\documentclass[11pt,a4paper]{article}
\usepackage{tikz}
\usepackage{subfigure}
\usepackage{enumerate}
\usepackage[english]{babel}
\usetikzlibrary{arrows, graphs,shapes}
\usepackage{amsfonts,amssymb,amsmath,latexsym,amsthm}
\usepackage{multirow}
\usepackage{bbm}
\usepackage{thm-restate}
\usepackage{algorithm}     
\usepackage{algorithmic}   
\usepackage{float}
\usepackage[T1]{fontenc}
\usepackage[utf8]{inputenc}

\usepackage[noabbrev,capitalize,nameinlink]{cleveref}

\newcommand \FF{\mathcal F}
\newcommand \GG{\mathcal G}
\newcommand \TT{\mathcal T}
\newcommand \kd{\operatorname{def}}
\newcommand \dist{\operatorname{dist}}

\newcommand \flev{\operatorname{flev}}
\DeclareFontFamily{U}{mathx}{}
\DeclareFontShape{U}{mathx}{m}{n}{<-> mathx10}{}
\DeclareSymbolFont{mathx}{U}{mathx}{m}{n}
\DeclareMathAccent{\fmz}{0}{mathx}{"71}

\newtheorem{thm}{Theorem}[section]
\newtheorem{cor}[thm]{Corollary}
\newtheorem{lem}[thm]{Lemma}
\newtheorem{prop}[thm]{Proposition}
\newtheorem{prob}[thm]{Problem}
\newtheorem{obs}[thm]{Observation}

\newtheorem{claim}{Claim}

\begin{document}

\title{Induced subgraphs of graphs with large deficiency }
\author{Jin Sun$^{a}$, \quad Xinmin Hou$^{a,b}$\\
\small $^{a}$ School of Mathematical Sciences\\
\small University of Science and Technology of China, Hefei 230026, Anhui, China\\
\small$^b$ Hefei National Laboratory\\
\small University of Science and Technology of China,  Hefei 230088, Anhui, China\\
\small Email: jinsun@mail.ustc.edu.cn (J. Sun), xmhou@ustc.edu.cn (X. Hou)
}
\date{}
\maketitle
\begin{abstract}
The deficiency of a graph $G$, denoted by $\kd(G)$, is the number of vertices not saturated by a maximum matching.
A bone $B_i$ is the tree obtained by attaching two pendent edges to each of the end vertices of a path $P_{i}$.  The local independence number of $G$, denoted by $\alpha_l(G)$, is defines as the maximum integer $t$ such that $G$ contains an induced star $K_{1,t}$.
Motivated by the seminal works of Scott and Seymour~(2016), Chudnovsky et al. (2017, 2020) on finding special types of holes in graphs with large chromatic number and bounded clique number, we establish an analog result by finding special types of bones in graphs with large deficiency and bounded local independence number. 
Fujita et al. (2006) proved  that $\kd(G)\le n-2$ if $G$ is a connected graph with $\alpha_l(G)<n$ and containing no bones.
We further establish exact extremal deficiency bounds for connected graphs with bounded local independence number that exclude specific bone configurations.
An algorithm that constructs large matchings and establishes an upper bound on the deficiency is also provided.
\end{abstract}
\section{Introduction}
In this paper, all graphs under consideration are finite and simple. Given a pair of graph parameters $\rho$ and $\sigma$, a class of graphs $\FF$ is called $(\rho,\sigma)$-bounded, if there exists a function $f$ so that each graph $G\in \FF$ satisfies $\rho(G)\le f(\sigma(G))$. The concept of  $(\rho,\sigma)$-boundedness has been extensively studied in the literature for various parameter pairs.
To summarize these results, we first introduce relevant graph parameters. 

For a graph $G$, let $\delta(G)$, $\omega(G)$, $\alpha(G)$, and $\chi(G)$ denote the minimum degree, clique number, independence number, and chromatic number of $G$,  respectively. 
We denote the {\it local independence number} of $G$ by $\alpha_l(G)$, which is defined as the largest integer $t$ such that $G$ contains an induced subgraph isomorphic to $K_{1,t}$.
We denote the {\it biclique number} of $G$ by $\tau(G)$, which is defined as the largest integer $t$ such that $G$ contains a subgraph (not necessarily induced) isomorphic to $K_{t,t}$. A {\it tree decomposition} of graph $G$ is a pair $\mathcal T=(T,\{B_t:t\in V(T)\})$, where $T$ is a tree and each $t\in V(T)$ is associated with a vertex subset $B_t\subseteq V(G)$. These subsets $B_t$ satisfy the following conditions:
\begin{itemize}
	\item For any $v\in V(G)$, there exists some $t\in V(T)$ such that $v\in B_t$,
	\item For any $uv\in E(G)$, there exists some $t\in V(T)$ such that $\{u,v\}\subseteq B_t$,
	\item For any $v\in V(G)$, the set $T_v=\{t\in V(T):v\in B_t\}$ induces a connected subtree of $T$.
\end{itemize}
The {\it width} of $\TT$, denoted by w$(\TT)$, and the {\it independence number} of $\TT$, denoted by $\alpha(\TT)$, is defined as follows respectively:
\begin{align*}
\text{w}(\TT):&= \max_{t\in V(T)} |B_t|-1,\\
\alpha(\TT):&= \max_{t\in V(T)} \alpha(G[B_t]).	
\end{align*}
An induced matching $M$ is said to {\it touch} a set $X\subseteq V(G)$ if every edge in $M$ has an endpoint in $X$. 
The {\it induced matching treewidth} of $\TT$, denoted by $\mu(\TT)$, is defined as the maximum, over 
all $t\in V(T)$, of the largest size of an induced matching touching $B_t$. 
The minimum values of w$(\TT)$, $\alpha(\TT)$ and $\mu(\TT)$ over all tree decomposition $\TT$ of $G$ are called the {\it treewidth}, the {\it tree independence number} and the {\it induced matching treewidth} of $G$, and are denoted by tw$(G)$, tree-$\alpha(G)$ and tree-$\mu(G)$ respectively. Clearly 
tree-$\mu(G)\le $ tree-$\alpha(G)\le $ tw$(G)+1$ 
since one can verify that  $\mu(\TT)\le \alpha(\TT)\le \text{w} (\TT)+1$ for any tree decomposition $\TT$ of $G$. 

In the following, we list some known results of $(\rho,\sigma)$-boundedness for various parameter pairs.
\begin{enumerate}[(1)]
 \item $(\delta, \chi)$-boundedness by Gy\'arf\'as, and Zaker~\cite{gz}.
 \item $(\text{tree-}\alpha,\text{tree-}\mu)$-boundedness by Abrishami, Bria\'nski, Czy\'zewska, McCarty, Milani\v c,  Rzą\'zewski, and Walczak~\cite{abcmmrw}.
 \item $(\text{tw},\omega)$-boundedness by Dallard, Milani\v c, and \v Storgel~\cite{dms}.
 \item $(\text{tree-}\alpha,\alpha_l)$-boundedness by Dallard, Krnc, Kwon, Milani\v c, Munaro, \v Storgel, and Wiederrecht~\cite{dkkmmsw}.
 \item $(\delta,\tau)$-boundedness, one can refer to the survey by Du, and McCarty~\cite{dm}.
 \item $(\chi,\omega)$-boundedness, one can refer to the surveys by Schiermeyer, Randerath \cite{sr}, and Scott, Seymour~\cite{ss20}.
\end{enumerate}	

Specially, the study of $(\chi,\omega)$-boundedness attracts long-standing interest. An {\em (odd) hole} in $G$ means an induced cycle of (odd) length at least four. 
In 1987, Gy\'arf\'as~\cite{g} proposed three conjectures to find holes in graphs with large chromatic number and bounded clique number, which are settled in the following.

\begin{thm}\label{long odd hole}
 (1) (Scott and Seymour~\cite{ss}) For all $k\ge 0$, there exists constant $c\ge 0$ such that for every graph $G$, if $\omega(G)\le k$ and $\chi(G) > c$ then $G$ has an odd hole.

 (2) (Chudnovsky, Scott, and Seymour~\cite{css}) For all $k,\ell \ge 0$, there exists constant $c\ge 0$ such that for every graph $G$, if $\omega(G)\le k$ and $\chi(G) > c$ then $G$ has a hole of length at least $\ell$.	

(3) (Chudnovsky, Scott, Seymour, and Spirkl~\cite{csss}) For all $k,\ell \ge 0$, there exists constant $c\ge 0$ such that for every graph $G$, if $\omega(G)\le k$ and $\chi(G) > c$ then $G$ has an odd hole of length at least $\ell$.	
\end{thm}

It is clear that (3) is stronger than (1) and (2). The {\it deficiency} of $G$, denoted by $\kd(G)$, is the number of vertices not saturated by a maximum matching. We will study $(\kd,\alpha_l)$-boundedness and give results analogous to \cref{long odd hole}.

For a given class of graphs $\FF$, a graph $G$ is called $\FF$-free if $G$ contains no induced copy from $\FF$. Graph $G$ is called {\it claw-free} if $G$ is $K_{1,3}$-free, i.e., $\alpha_l(G)<3$. 
A classical result of (near-)perfect matching in claw-free graphs due to Las Vergnas \cite{l}, Sumner \cite{s},  J\"unger, Pulleyblank and Reinelt~\cite{jpr} can be described as follows.
 
\begin{thm}[\cite{jpr,l,s}]\label{sumner}
Let $G$ be a connected graph with $\alpha_l(G)<3$. Then $\kd(G)\le 1$. 
 \end{thm}

Let $P_i$ be a path with $i$ vertices. Let $i\ge 2$.
Corresponding to a hole, {\em a bone} $B_i$ is a tree formed by attaching two pendent edges to each end vertex of a path $P_{i}$. When $i$ is odd (even), $B_i$ is called an odd (even) bone. We call a set $A\subseteq \{2,3,4,5\cdots\}$  an {\it admitting set}. 
Let $\GG(A)$ denote the family of connected graphs in  which all bones belong to $\{B_i : i\in A\}$. 
Let $A_{\ge m}=\{a\in A:a\ge m\}$.

When all bones are forbidden in graphs, Fujita, Kawarabayashi, Lucchesi, Ota, Plummer and Saito~\cite{fklops} extended \cref{sumner} to  bone-free connected graphs, showing that:
\begin{thm}[\cite{fklops}]\label{bone}
Given integer $n>3$, let $G\in \GG(\emptyset)$ satisfies $\alpha_l(G)<n$. Then $\kd(G)\le n-2$.	
\end{thm}

In this paper, we further strengthen \cref{bone} to the family of connected graphs with bounded local independence number which admit a special set of bones. 
\begin{restatable}{thm}{main}
\label{main}
Let $m\ge 3$ be an odd integer and let $A$ be a set of odd integers such that for any $p,q\in A$ with $p,q\ge m$, neither $p+q+1$ nor $p+q-1$ belongs to $A$. For any integer $n>3$, every graph $G$ in $\GG(A)$ with $\alpha_l(G)<n$ satisfies $\kd(G)\le m(n-3)(n-2)^{ \frac{m-3}{2} }+1$.
Furthermore, when $m=3$, the bound improves to $\kd(G)\le 2n-5$.
\end{restatable}

\noindent{\bf Remark:} The upper bound $\kd(G)\le m(n-3)(n-2)^{ \frac{m-3}{2} }+1$ is asymptotically optimal with respect to the leading item $(n-3)(n-2)^{ \frac{m-3}{2} }$, and when $m=3$, the upper bound $\kd(G)\le 2n-5$ is tight. The tightness of these upper bounds are shown in the following constructions.

\begin{figure}
\centering
\begin{tikzpicture}[scale=0.9]
\begin{scope}
\fill (0,-1) circle (1.8pt);\fill (0,1) circle (1.8pt);\fill (0,0) circle (1.8pt);\fill (1,0) circle (1.8pt);
\node at (0,0.6) {$\vdots$};
\draw (-0.5,-1.5) rectangle (0.5,1.5);\node at (0,-0.5) {$E_n$};
			
\fill (5,-1) circle (1.8pt);\fill (5,1) circle (1.8pt);\fill (5,0) circle (1.8pt);\fill (4,0) circle (1.8pt);
\node at (5,0.6) {$\vdots$};
			
\fill (1.5,0)node[below]{$2$} circle (1.8pt);
			
\draw (1,0)edge (0,1) edge(0,0) --(0,-1);\draw(1,0) node[below]{$1$} -- (2.3,0);
\draw(3.2,0) -- (4,0)node[below]{$p$}; \draw (4,0) edge (5,-1) edge(5,1) --(5,0);
\node at (2.8,0) {$\dots$};
\node at (2.5,-2.5) {$B\!S_n^p$};	
\end{scope}

\begin{scope}[shift={(8,-1.5)}]
\fill (-0.5,1.5) circle (1.8pt);
\draw (-0.5,1.5) edge(0.5,0.5) edge(0.5,1.5) --(0.5,2.5);
\draw (0,0) rectangle (1,3);
\fill (0.5,0.5) circle (1.8pt);
\node at (0.5,2.1) {$\vdots$}; 
\fill (0.5,1.5) circle (1.8pt); 
\fill (0.5,2.5) circle (1.8pt)node[left]{$1$}; 
\node at (0.5,1) {$E_n$};

\fill (4,0.5) circle (1.8pt);
\draw (0.5,0.5)--(4,0.5); 
\fill (0.5,1.5) circle (1.8pt); 
\fill (4,1.5) circle (1.8pt);
\draw (0.5,1.5)--(4,1.5);

\fill (0.5,2.5) circle (1.8pt);
\fill (1.5,2.5) circle (1.8pt)node[below]{$2$};
\fill (4,2.5) circle (1.8pt)node[below,red]{$p\!-\!1$};

\draw (6,2.5)node [red] {$E_{n-1}$} circle (1 and 0.3);
\draw (6,1.5) circle (1 and 0.3);
\draw (6,0.5) circle (1 and 0.3);
\draw[blue, ultra thick] (4,2.5)--(5,2.5);
\draw[blue, ultra thick] (4,1.5)--(5,1.5);
\draw[blue, ultra thick] (4,0.5)--(5,0.5);

\draw (0.5,2.5)--(4,2.5);
\node at (2.5,2.1) {$\vdots$};
\node at (3,-1) {$S_n^p$};		
\end{scope}

\begin{scope}[shift={(0,-5)}]
\fill (0,-1) circle (1.8pt);\fill (0,1) circle (1.8pt);\fill (0,0) circle (1.8pt);\fill (1,0) circle (1.8pt);
\node at (0,0.6) {$\vdots$};
\draw (-0.5,-1.5) rectangle (0.5,1.5);\node at (0,-0.5) {$E_n$};
			
\fill (5,-1) circle (1.8pt);\fill (5,1) circle (1.8pt);\fill (5,0) circle (1.8pt);\fill (4,0) circle (1.8pt);
\node at (5,0.6) {$\vdots$};
			
\fill (1.5,0)node[below]{$x_2$} circle (1.8pt);
			
\draw (1,0)edge (0,1) edge(0,0) --(0,-1);\draw(1,0) node[below]{$x_1$} -- (2.3,0);
\draw(3.2,0) -- (4,0)node[below]{$x_p$}; \draw (4,0) edge (5,-1) edge(5,1) --(5,0);
\node at (2.8,0) {$\dots$}; 

\fill (10,-1) circle (1.8pt);\fill (10,0) circle (1.8pt);\fill (10,1) circle (1.8pt);\fill (9,0) circle (1.8pt);
\node at (10,0.6) {$\vdots$};
			
\fill (6.5,0)node[below]{$y_2$} circle (1.8pt);\fill (6,0)node[below]{$y_1$} circle (1.8pt);
			
\draw (6,0)edge (5,1) edge(5,0) --(5,-1);\draw(6,0) -- (7.3,0);
\draw(8.2,0) -- (9,0)node[below]{$y_p$}; \draw (9,0) edge (10,-1) edge(10,1) --(10,0);
\node at (7.8,0) {$\dots$};    
\node at (5,-2) {$T_n^p$};
\end{scope}

\begin{scope}[shift={(-0.5,-12)}]
\fill[gray!35] (0,0) rectangle (1,3);
\fill (0.5,0.5) circle (1.8pt);
\node at (0.5,2.1) {$\vdots$}; 
\fill (0.5,1.5) circle (1.8pt); 
\fill (0.5,2.5) circle (1.8pt)node[left]{$1$}; 
\node at (0.5,1) {$K_m$};

\fill (4,0.5) circle (1.8pt);
\draw (0.5,0.5)--(4,0.5); 
\fill (0.5,1.5) circle (1.8pt); 
\fill (4,1.5) circle (1.8pt);
\draw (0.5,1.5)--(4,1.5);

\fill (0.5,2.5) circle (1.8pt);
\fill (1.5,2.5) circle (1.8pt)node[below]{$2$};
\fill (4,2.5) circle (1.8pt)node[below]{$p$};

\draw (6,2.5)node {$E_n$} circle (1 and 0.3);
\draw (6,1.5) circle (1 and 0.3);
\draw (6,0.5) circle (1 and 0.3);
\draw[blue, ultra thick] (4,2.5)--(5,2.5);
\draw[blue, ultra thick] (4,1.5)--(5,1.5);
\draw[blue, ultra thick] (4,0.5)--(5,0.5);

\draw (0.5,2.5)--(4,2.5);
\node at (2.5,2.1) {$\vdots$};
\node at (3.5,-1) {$E_{m,n}^p$};	
\end{scope}

\begin{scope}[shift={(8,-12)}]
\fill[gray!30] (0,0) rectangle (1,3);
\fill (0.5,0.5) circle (1.8pt);
\node at (0.5,2.1) {$\vdots$}; 
\fill (0.5,1.5) circle (1.8pt); 
\fill (0.5,2.5) circle (1.8pt)node[left]{$1$}; 
\node at (0.5,1) {$K_m$};

\fill (4,0.5) circle (1.8pt);
\draw (0.5,0.5)--(4,0.5); 
\fill (0.5,1.5) circle (1.8pt); 
\fill (4,1.5) circle (1.8pt);
\draw (0.5,1.5)--(4,1.5);

\fill (0.5,2.5) circle (1.8pt);
\fill (1.5,2.5) circle (1.8pt)node[below]{$2$};
\fill[red] (1.5,3.5) circle (1.8pt)node[right]{$v$};
\draw (1.5,3.5) edge(0.5,2.5)--(1.5,2.5);
\fill (4,2.5) circle (1.8pt)node[below]{$p$};

\draw (6,2.5)node {$E_n$} circle (1 and 0.3);
\draw (6,1.5) circle (1 and 0.3);
\draw (6,0.5) circle (1 and 0.3);
\draw[blue, ultra thick] (4,2.5)--(5,2.5);
\draw[blue, ultra thick] (4,1.5)--(5,1.5);
\draw[blue, ultra thick] (4,0.5)--(5,0.5);

\draw (0.5,2.5)--(4,2.5);
\node at (2.5,2.1) {$\vdots$};
\node at (3.5,-1) {$E_{m,n}^{p+}$};	
\end{scope}
\end{tikzpicture}	
\caption{Graphs $B\!S_n^p, S_n^p,T_n^p,E_{m,n}^p$ and $E_{m,n}^{p+}$. All edges exist between the two end of the thick blue line. The red vertex $v$ is the vertex added to $E_{m,n}^{p+}$ from $E_{m,n}^p$.}
  \label{fig 1}
\end{figure}
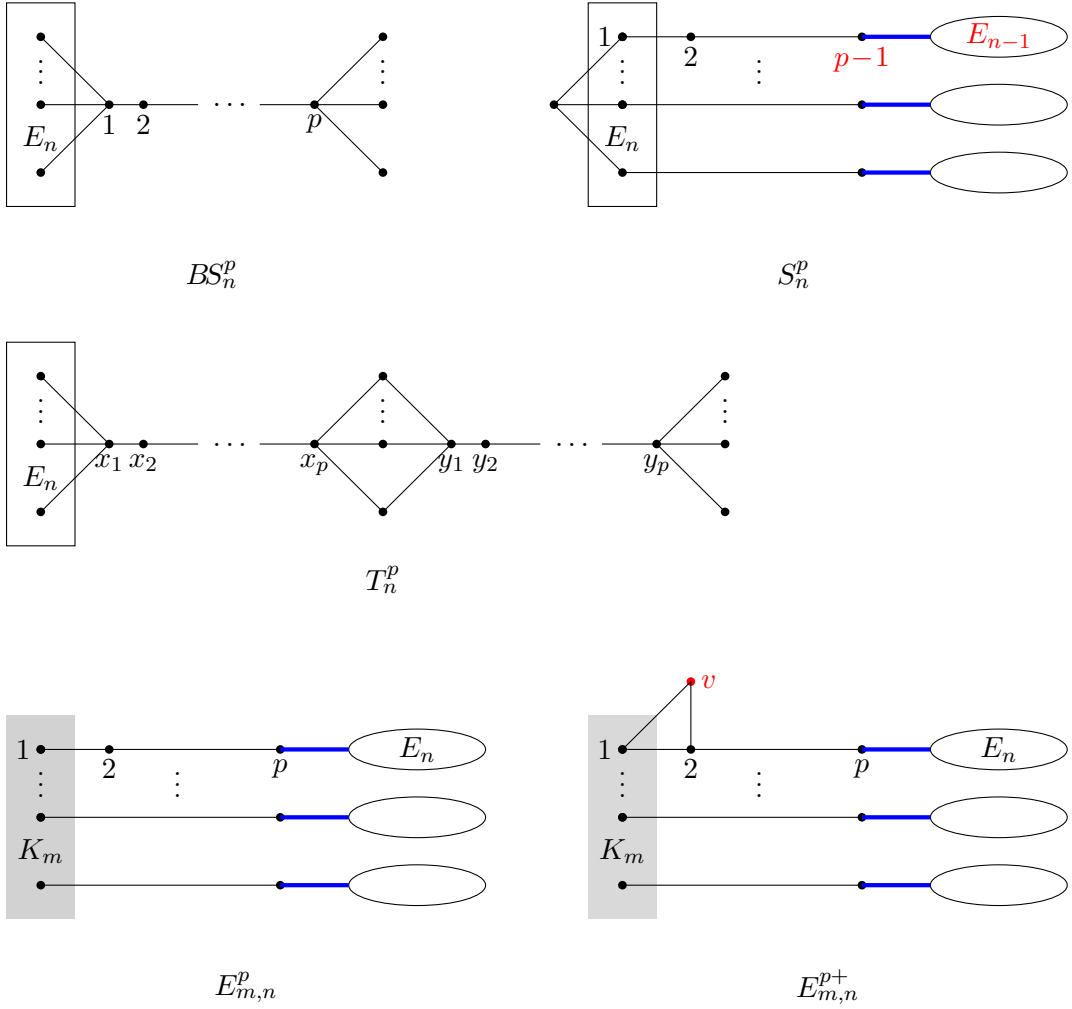

\noindent{\bf Construction A}: We call an edge $xy$ a {\em pendent edge} of vertex $x$ if $N(y)=\{x\}$. Let $n, p$ be positive integers. Let $D_n^p$ be the tree formed by attaching $n$ pendent edges at an end vertex of a path $P_p$. We regard the other end of $P_p$ as the end of $D_n^p$. Specially, $D_2^1$ is a path $P_3$ with the middle vertex as its end. Attaching $D_n^p$ at a vertex $v$ of a graph $G$ means that we identify the end of $D_n^p$ with $v$. 
Several extremal graphs are constructed as follows (\cref{fig 1}).

\begin{itemize}
\item[$B\!S_n^p$:] the tree formed by attaching $n$ pendent edges to each end vertex of a path $P_{p}$. Actually, a bone $B_p$ is a $B\!S_2^p$. 
	\item[$S_n^p:$] the tree formed by attaching $n$ pairwise disjoint copies of $D_{n-1}^p$ at a common vertex.
	\item[$T_n^p:$] the graph obtained by attaching one copy of $D_n^p$ at each of the two vertices of degree $n$ in $K_{2,n}$.
 \item[$E_{m,n}^p:$] the graph formed by attaching a copy of $D_n^p$ at each vertex of the complete graph $K_m$.
\item[$E_{m,n}^{p+}$:] the graph obtained from $E_{m,n}^p$ by adding a new vertex adjacent to two consecutive vertices at the end of one of the $D_n^p$ copy.
\end{itemize}

When $p$ is odd, $B\!S_{n-2}^p$ is $K_{1,n}$-free with admitting set $A=\{p\}$. Clearly, $x+y\pm1\notin A$ for $x,y\in A$ with $x,y\ge 3$. It can be checked directly that $\kd(BS_{n-2}^p)=2n-5$. This implies that the upper bound $\kd(G)\le 2n-5$ is tight when $m=3$.

\bigskip 
\noindent{\bf Construction B}: Let $T_{3,n}=K_{1,n-1}$. Given an odd integer $m\ge 5$ and an integer $n>3$, define $T_{m,n}$ as a rooted tree with vertex degrees restricted to $\{1,2,n-1\}$, structured as follows: The vertices are partitioned into levels $L_0, L_1, L_2, \dots, L_{m-2}$ based on their distances from the root. The degree arrangements are then constructed as:   
\begin{align*}
V_1(T_{m,n}):&=\{v\in V(T):\deg(v)=1\}=L_{m-2}, \\
V_2(T_{m,n}):&=\{v\in V(T):\deg(v)=2\}=L_1\cup L_3 \cup \cdots \cup L_{m-4},\\
V_{n-1}(T_{m,n}):&=\{v\in V(T):\deg(v)=n-1\}=L_0\cup L_2 \cup \cdots \cup L_{m-3}.
\end{align*}
Thus $V(T_{m,n})=V_1(T_{m,n})\cup V_2(T_{m,n})\cup V_{n-1}(T_{m,n})$. 

$\alpha_l(T_{m,n})<n$ as the maximum degree of $T_{m,n}$ is $n-1$. We now consider the admitting set of $T_{m,n}$. In order to find a bone in tree $T_{m,n}$, we need to find vertices $u_1$ and $u_2$ satisfying $\deg(u_i)\ge 3$ for $i=1,2$.  Hence  $u_i$ must be in $V_{n-1}(T_{m,n})=L_0\cup L_2 \cup \cdots \cup L_{m-3}$ for $i=1,2$. Since there is exactly one path connecting $u_1$ and $u_2$ in $T_{m,n}$, the possible bone obtained from this path must be $B_{\dist(u_1,u_2)+1}$. Therefore, $T_{m,n}$ has admitting set $A=\{3,\cdots , 2m-5\}$. Thus $T_{m,n}\in\GG(A)$.
Clearly, $p+q\pm 1\notin A$ for any $p,q\in  A_{\ge m}$. 

We can calculate deficiency of graphs with pendent edges effectively according to the following observation.
\begin{obs}\label{pendent edge}
Given a graph $G$, let $uv$ be a pendent edge of vertex $u$. Then $\kd(G)=\kd(G-\{u,v\})$.
\end{obs}
\begin{proof}
Let $M$ be a maximum matching of $G$. Then $u\in V(M)$; otherwise, $v\not \in V(M)$ neither. Thus $M\cup \{uv\}$ is a larger matching, a contradiction. If $uw\in M$ for some $w$ distinct from $v$, then $M'=M-\{uw\}\cup \{uv\}$ is still a maximum matching. Hence $G$ contains a maximum matching $M'$ satisfying $uv\in M'$. This leads to $\kd(G)=\kd(G-\{u,v\})$.
\end{proof}
\begin{cor}\label{Tree}
$\kd(T_{m,n})=(n-1)(n-2)^{\frac{m-3}{2}}-1$.  	
\end{cor}
\begin{proof}
In tree $T_{m,n}$, we can always find a pendent edge $uu^+$ for each vertex $u\in L_{m-3}$ with $u^+\in L_{m-2}$. 
Note that $|L_i|=(n-1)(n-2)^{\frac{i-2}2}$ for even $i\ge 2$. Thus
$|L_{m-3}|=(n-1)(n-2)^{\frac{m-5}{2}}$. Hence the resulting graph $T_{m,n}-\{u,u^+:u\in L_{m-3}\}$ is a disjoint union of $T_{m-2,n}$ and $(n-3)|L_{m-3}|=(n-3)(n-1)(n-2)^{\frac{m-5}{2}}$ isolated vertices. 
By \cref{pendent edge}, 
\[\kd(T_{m,n})=\kd(T_{m,n}-\{u,u^+:u\in L_{m-3}\})=\kd(T_{m-2,n})+(n-3)(n-1)(n-2)^{\frac{m-5}{2}}.\] 
By induction, we have  \[\kd(T_{m,n})=\kd(T_{3,n})+(n-1)(n-3)\sum_{i=0}^{\frac{m-5}{2}}(n-2)^i=(n-1)(n-2)^{\frac{m-3}{2}}-1.\]	
\end{proof}

By \cref{Tree}, the upper bound $\kd(G)\le m(n-3)(n-2)^{ \frac{m-3}{2} }+1$ is asymptotically optimal with respect to the leading item $(n-3)(n-2)^{ \frac{m-3}{2} }$. 
It is interesting to determine the exact upper bound of $\kd(G)$. 

\begin{prob}
Given odd integer $m\ge 5$ and integer $n>3$, determine the exact upper bound of $\kd(G)$ for graph $G$ with $\alpha_l(G)<n$ and admitting set $A$ consisting of odd integers such that for any $p,q\in A$ with $p,q\ge m$, $p+q\pm 1\not \in A$.
\end{prob}

Several extremal results concerning $(\chi,\omega)$-boundedness of graphs that forbid certain  types of  holes have been established, as seen in  \cite{crst,cs,wxx}. A connected graph $G$ is called {\it deficiency-critical} if for any proper connected induced subgraph $G'$ of $G$, we have $\kd(G')<\kd(G)$. For example, all odd bones are deficiency-critical while even bones not. 
The following theorem demonstrates the necessity of the condition `for all $ p,q \in A_{\ge m}$, $p+q\pm 1\not \in A$' in \cref{main}.

\begin{restatable}{thm}{twooddbone}
\label{two odd bone}
Given an integer $n>3$ and an odd integer $p\ge 3$, let $G\in \GG(\{p,q\})$ satisfies $\alpha_l(G)<n$.
\begin{itemize}
\item If $q=2p+1$, then $\kd(G)\le 3n-8$. Moreover, the equality holds if and only if $G\cong T_{n-2}^p$ in the context of deficiency-critical graphs.	
\item If $q=2p-1$, then $\kd(G)\le n^2-3n+1$. Moreover,  the equality holds if and only if $G\cong S_{n-1}^p$ in the context of deficiency-critical graphs.	
\end{itemize}
\end{restatable}

When the admitting set $A$ contains even integers, the situation will be different. Graph $E_{m,2}^p$ has admitting set $\{2p\}$ and is $K_{1,4}$-free. However, $\kd(E_{m,2}^{p})\ge m$ can be arbitrarily large along with the increasing of $m$. Therefore, in order to bound the deficiency, we need additional restriction on the clique number $\omega(G)$.  

\begin{figure}
	\centering
	\begin{tikzpicture}[scale=1.4]

\node[fill=blue, regular polygon, regular polygon sides=3, inner sep=2pt] at (0,0) {};
\fill (-1,1) circle (1.5pt);\fill (0,1) circle (1.5pt);\fill (1,1) circle (1.5pt);
\draw (0,0) edge(-1,1) edge(0,1) --(1,1);
\fill (-2,2) circle (1.5pt);\fill (0,2) circle (1.5pt);\fill (2,2) circle (1.5pt);
\draw[dashed] (-1,1)--(-2,2) (0,1)--(0,2) (1,1)--(2,2);
\node[fill=blue, regular polygon, regular polygon sides=3, inner sep=2pt] at (-2,3) {};
\node[fill=blue, regular polygon, regular polygon sides=3, inner sep=2pt] at (0,3) {};
\node[fill=blue, regular polygon, regular polygon sides=3, inner sep=2pt] at (2,3) {};
\draw (-2,2)--(-2,3) edge(-2.5,4)--(-1.5,4) (0,2)--(0,3)edge(-0.5,4)--(0.5,4) (2,2)--(2,3)edge (1.5,4)--(2.5,4);
\fill (-2.5,4) circle (1.5pt);\fill (-1.5,4) circle (1.5pt);\fill (-0.5,4) circle (1.5pt);\fill (0.5,4) circle (1.5pt);\fill (1.5,4) circle (1.5pt);\fill (2.5,4) circle (1.5pt);
\draw[dashed] (-2.5,4)--(-2.5,5)edge(-2.2,6)--(-2.8,6) (-1.5,4)--(-1.5,5)edge(-1.8,6)--(-1.2,6) (-0.5,4)--(-0.5,5)edge(-0.2,6)--(-0.8,6) (0.5,4)--(0.5,5)edge(0.2,6)--(0.8,6) (1.5,4)--(1.5,5)edge(1.2,6)--(1.8,6) (2.5,4)--(2.5,5)edge(2.2,6)--(2.8,6);
\node[fill=blue, regular polygon, regular polygon sides=3, inner sep=2pt] at (-2.5,5) {};
\node[fill=blue, regular polygon, regular polygon sides=3, inner sep=2pt] at (-1.5,5) {};
\node[fill=blue, regular polygon, regular polygon sides=3, inner sep=2pt] at (-0.5,5) {};
\node[fill=blue, regular polygon, regular polygon sides=3, inner sep=2pt] at (0.5,5) {};
\node[fill=blue, regular polygon, regular polygon sides=3, inner sep=2pt] at (1.5,5) {};
\node[fill=blue, regular polygon, regular polygon sides=3, inner sep=2pt] at (2.5,5) {};

\fill (-2.25,6) circle (1.5pt);\fill (-2.75,6) circle (1.5pt);\fill (-1.75,6) circle (1.5pt);\fill (-1.25,6) circle (1.5pt);\fill (-0.75,6) circle (1.5pt);\fill (-0.25,6) circle (1.5pt);\fill (0.25,6) circle (1.5pt); \fill (0.75,6) circle (1.5pt);\fill (1.25,6) circle (1.5pt);\fill (1.75,6) circle (1.5pt);\fill (2.25,6) circle (1.5pt);\fill (2.75,6) circle (1.5pt);

\fill (-2.35,6.25) circle (1.5pt);\fill (-2.85,6.25) circle (1.5pt);\fill (-1.85,6.25) circle (1.5pt);\fill (-1.35,6.25) circle (1.5pt);\fill (-0.85,6.25) circle (1.5pt);\fill (-0.35,6.25) circle (1.5pt);\fill (0.15,6.25) circle (1.5pt); \fill (0.65,6.25) circle (1.5pt);\fill (1.15,6.25) circle (1.5pt);\fill (1.65,6.25) circle (1.5pt);\fill (2.15,6.25) circle (1.5pt);\fill (2.65,6.25) circle (1.5pt);
\fill (-2.15,6.25) circle (1.5pt);\fill (-2.65,6.25) circle (1.5pt);\fill (-1.65,6.25) circle (1.5pt);\fill (-1.15,6.25) circle (1.5pt);\fill (-0.65,6.25) circle (1.5pt);\fill (-0.15,6.25) circle (1.5pt);\fill (0.35,6.25) circle (1.5pt); \fill (0.85,6.25) circle (1.5pt);\fill (1.35,6.25) circle (1.5pt);\fill (1.85,6.25) circle (1.5pt);\fill (2.35,6.25) circle (1.5pt);\fill (2.85,6.25) circle (1.5pt);

\draw (-2.25,6)edge +(0.1,0.25) -- +(-0.1,0.25) (-2.75,6)edge +(0.1,0.25) -- +(-0.1,0.25) (-1.75,6)edge +(0.1,0.25) -- +(-0.1,0.25) (-1.25,6)edge +(0.1,0.25) -- +(-0.1,0.25) (-0.75,6)edge +(0.1,0.25) -- +(-0.1,0.25) (-0.25,6)edge +(0.1,0.25) -- +(-0.1,0.25) (0.25,6)edge +(0.1,0.25) -- +(-0.1,0.25) (0.75,6)edge +(0.1,0.25) -- +(-0.1,0.25) (1.25,6)edge +(0.1,0.25) -- +(-0.1,0.25) (1.75,6)edge +(0.1,0.25) -- +(-0.1,0.25) (2.25,6)edge +(0.1,0.25) -- +(-0.1,0.25) (2.75,6)edge +(0.1,0.25) -- +(-0.1,0.25);
\node at (-4,0){$L_0$}; \node at (-4,1){$L_1$}; \node at (-4,2){$L_{a_1-1}$};\node at (-4,3){$L_{a_1}$};\node at (-4,4){$L_{a_1+1}$};\node at (-4,5){$L_{a_2}$};\node at (-4,6){$L_{a_3}$};
\end{tikzpicture}
\caption{$F(a_1,a_2,a_3)$: the dashed line represents path,
the blue triangle in $L_0\cup L_{a_1}\cup L_{a_2}$ means the $K_3$ got by $Y$-$\Delta$ operation.
}
	\label{fig 2}
\end{figure}
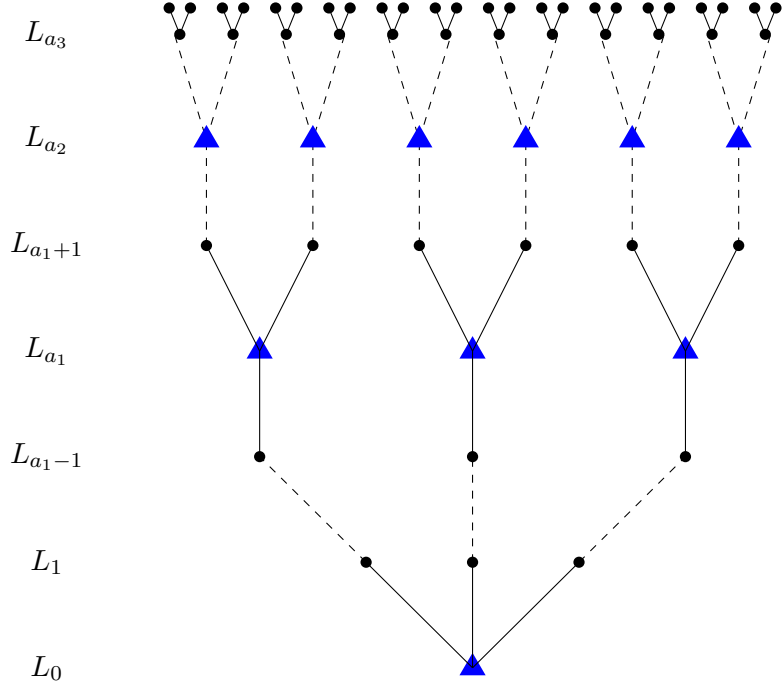

\noindent{\bf Construction C}: For a graph $G$ and a vertex $v\in V(G)$ with neighborhood $N(v)=\{a,b,c\}$, a {\em $Y-\Delta$ operation} at $v$ consists of replacing $v$ with three new vertices $\{x,y,z\}$ and replacing $\{va,vb,vc\}$ with six new edges  $\{xy,yz,zx,ax,by,cz\}$.
Given positive integers $a_1< a_2<\cdots <a_i$, define a rooted tree $T=T(a_1,\cdots, a_i)$ with vertex degrees restricted to $\{1,2,3\}$, constructed as follows: The vertices are partitioned into levels $L_0, L_1, \dots, L_{a_i}$ based on their distances from the root. The degree arrangements are constructed as: 
\begin{align*}
V_1(T):&=\{v\in V(T):\deg(v)=1\}=L_{a_i}, \\
V_2(T):&=\{v\in V(T):\deg(v)=2\}=\bigcup_{j=1}^{a_1-1}L_j\cup \bigcup_{j=a_1+1}^{a_2-1}L_j\cup \cdots \cup \bigcup_{j=a_{i\!-\!1}+1}^{a_i-1}L_j,\\
V_3(T):&=\{v\in V(T):\deg(v)=3\}=L_0\cup L_{a_1} \cup \cdots \cup L_{a_{i-1}}.
\end{align*}
Thus $V(T)=V_1(T)\cup V_2(T)\cup V_3(T)$.
Define $F(a_1,\ldots ,a_r)$ to be the graph obtained from $T(a_1,\ldots, a_r)$ by $Y-\Delta$ operations at each vertex of $V_3(T)$ and attaching a $D_2^1$ at each vertex of $V_1(T)$. (\cref{fig 2})

It can be checked directly that $F(a_1,\dots, a_r)$ is $(K_4,K_{1,4})$-free and contains only even bones, while  $\kd(F(a_1,\dots, a_r))\ge 3\cdot 2^{r-1}=|L_{a_i}|$ can be arbitrarily large along with the increasing of $r$. In order to find a bone in $F(a_1,\dots, a_r)$, we need to find vertices $u_1$ and $u_2$ satisfying $\alpha(N(u_i))\ge 3$ for $i=1,2$. Thus $u_i\in L_{a_i}$. The induced path in $F(a_1,\dots, a_r)$ connecting $u_1$ and $u_2$ is unique, which has length $\dist_{T}(u_1,u_2)+1$. Thus the possible bone found by $u_1$ and $u_2$ will be $B_{2a_i-2a_j+2}$ for some $j$ with $0\le j \le i-1$. (Let $a_0=0$ for convenience) Therefore, $F(a_1,\dots, a_r)$ has admitting set $A=\{2a_i-2a_{i-1}+2,2a_i-2a_{i-2}+2,\cdots , 2a_i-2a_1+2,2a_i+2 \}$. $A$ can be arbitrarily sparse by the choices of reasonably positive integers $a_1< a_2<\cdots <a_i$. 

Thus, we require that $G$ contains no $K_3$ or $G$ contains only one type of even bone in the following.
\begin{restatable}{thm}{onlyoneevenboneevenbone}
\label{only one even bone}
(1) Given integers $m,n> 3$ and an even integer $2p\ge 2$, if $G\in \GG(\{2p\})$ is $\{K_{1,n},K_m\}$-free, then $\kd(G)\le (m-1)(n-3)+1$. 

(2) Given an integer $n>3$, let $A$ be the set of all positive even integers.  If $G\in \GG(A)$ is $\{K_{1,n},K_3\}$-free, then $\kd(G)\le 2n-6$.
\end{restatable}

\noindent{\bf Remark}: (1) The upper bounds $\kd(G)\le (m-1)(n-3)+1$ is tight. If $(m-1)(p-1)$ is odd, then $E_{m-1,n-2}^p$ is a $\{K_m,K_{1,n}\}$-free graph with admitting set $\{2p\}$. We can calculate that $\kd(E_{m-1,n-2}^p)=(m-1)(n-3)+1$ by applying \cref{pendent edge} iteratively. On the other hand, if $(m-1)(p-1)$ is even, then $E_{m-1,n-2}^{p+}$ is a $\{K_m,K_{1,n}\}$-free graph with admitting set $\{2p\}$. At this time, it holds that $\kd(E_{m-1,n-2}^{p+})=(m-1)(n-3)+1$.

(2) The upper bound $\kd(G)\le 2n-6$ is tight as $B\!S_{n-2}^2$ is a $\{K_3,K_{1,n}\}$-free graph with admitting set $\{2\}$ and $\kd(B\!S_{n-2}^2)=2n-6$ attains this extremum.

The rest of this article is arranged as follows. In Section $2$, we introduce preliminary concepts, build an algorithm to construct matching and bound deficiency of graphs, present some preliminary lemmas. We prove Theorem~\ref{main} in Section 3 and, Theorems \ref{two odd bone} and \ref{only one even bone} in the last section.

\section{Preliminaries}
For positive integers $m$ and $n$ with $m<n$, let $[m,n]$ denote the set $\{m,m+1,\dots ,n-1,n\}$, and let $[n]$ denote $[1,n]$.
For a vertex $v$ in $V(G)$, let $N_G(v)=\{u\in V(G):uv\in E(G)\}$ be the open {\em neighborhood} of $v$ and let $N_G[v]=\{v\}\cup N_G(v)$ be the {\em closed neighborhood} of $v$. We will omit the subscript if there is no ambiguity. For a subset $X$ of $V(G)$, let $N(X)=\bigcup_{x\in X}N(x)-X$. Let $G[X]$ be the subgraph of $G$ induced by $X$. We call $X$ an independent set (or a stable set) if there is no edge in $G[X]$. For a path $P=x_1\cdots x_n$ of $G$, we use $x_i P x_j$ to denote the subpath from $x_i$ to $x_j$. The {\em distance} from a vertex $x$ to a vertex $y$ in a graph $G$, denoted by $\dist(x,y)$, is the length of a shortest $(x,y)$-path in $G$. 

Let $G$ be a connected graph and  $x$ be a vertex of $G$. We call the following partition of $V(G)$ as {\it levelling} of $G$ from $x$. For $i\ge 0$, let $L_i=\{u\in V(G):\dist(u,x)=i\}$. The maximum number $N$ such that $L_N\not =\emptyset$ is called the {\it level number}. Since $G$ is connected, we have partition $V(G)=\bigcup_{i=0}^{N}L_i$. If there is an edge $uv\in E(G)$, where $u\in L_i, v\in L_j$, then $|i-j|\le 1$ by the levelling construction. For any $u\in L_i$, we define $C(u)=N(u)\cap L_{i+1}$ as the {\em children} of $u$. Acquiescently, we always use $P_u=u_iu_{i-1}\cdots u_1u_0$ to denote a selected shortest path connecting $u$ and $x$, where $u_i=u$ and $u_0=x$. Let $v\in L_j$ be distinct from $u$. The {\em friendly level} of $u$ and $v$, denoted by $\flev(u,v)$, is the maximum integer $k$ such that there exists $P_u$ and $P_v$ satisfying $u_{k'}=v_{k'}$ for any $k'\le k$. It is well-defined since $u_0=v_0=x$ always. In this case, we call $P_u$ and $P_v$ are friendly at level $\flev(u,v)$. Thus $\{u_{k'}v_{k'+1},v_{k'}u_{k'+1}:k'>\flev(u,v)\}\cap E(G)=\emptyset$. For convenience, we abuse the notation even if $u_{k'}$ or $v_{k'}$ does not exist.

The following lemma will be used to construct matching from consecutive two levels of the levelling.
Let $G$ be a graph with vertex partition $V(G)=X\cup Y$. For any matching $M$ in $G$, define $G_M=G-V(M)$,
$Y_M=Y-V(M)$ (the set of vertices in $Y$ not saturated by $M$),
$X_M=\{x\in X-V(M) : N(x)\cap Y_M\neq \emptyset\}$ (the set of vertices in $X$ with neighbors in $Y_M$ and not saturated by $M$). Thus $Y_M\subseteq N(X_M)$ is equivalent to $Y_M\subseteq N(X-V(M))$.

\begin{lem}\label{two level}
If  $G$ is a graph with vertex partition $V(G)=X\cup Y$ satisfying $Y\subseteq N(X)$,
then there exists a matching $M$ in $G$ satisfying:
\begin{itemize}
	\item[(1)] $Y_M\subseteq N(X_M)$,
	\item[(2)] $Y_M$ is a stable set,
	\item[(3)] For any $x\in X_M$, there exist $y_1,y_2\in Y_M$, such that $N_{G_M}(y_1)=N_{G_M}(y_2)=\{x\}$,
	\item[(4)] For any $v\in V(M)\cap X$, at most one vertex $x\in X_M$ does not satisfy that, there exist $y_{1},y_{2}\in Y_M$ such that $N_{G_M\cup \{v\}}(y_1)=N_{G_M\cup \{v\} }(y_2)=\{x\}$, where $G_M\cup\{v\}=G[V(G_M)\cup \{v\}]$.
\end{itemize}
\end{lem}

\begin{proof}
Take a matching $M$ as large as possible satisfying $Y_M\subseteq N(X_M)$. We claim that (2), (3), and (4) hold.

Suppose to the contrary of $(2)$ that there exists an edge $uv$ in $Y_M$, then $M'=M\cup \{uv\}$ is a larger matching satisfying $Y_{M'}\subseteq N(X_{M'})$, a contradiction.

As for $(3)$, take $y\in N(x)\cap Y_{M}$  first for $x\in X_M$. Since $M'=M\cup \{xy\}$ will lead to $Y_{M'}\not \subseteq N(X_{M'})$, there must exist a vertex $y_1\in Y_{M}$ such that $N(y_1)\cap X_M=\{x\}$. 
Similarly, since $M'=M\cup \{xy_1\}$ will lead to $Y_{M'}\not \subseteq N(X_{M'})$, there exists another vertex $y_2\in Y_M$ such that $N(y_2)\cap X_M=\{x\}$.

Suppose to the contrary of $(4)$ that $u$ destroys $x_i\in X_M$ for $i=1,2$, i.e., at most one $y_i\in Y_M$ satisfies $N_{G_M\cup \{u\} }(y_i)=\{x_i\}$ for $i=1,2$. If there is no such vertex, let $y_i$ be an arbitrary vertex satisfying $N_{G_M }(y_i)=\{x_i\}$. Let $uv\in M$. Thus $M'=(M\setminus\{uv\})\cup \{x_1y_1,x_2y_2\}$ is a larger matching satisfying $Y_{M'}\subseteq N(X_{M'})$.
\end{proof}

The graph with a single vertex is deficiency-critical according to the definition of deficiency-critical graph. A graph $G$ is called nontrivial if $|V(G)|\ge 2$. 

\begin{cor}\label{snail horn}
Let $V(G)=\bigcup_{i=0}^{N}L_i$ be the levelling of a nontrivial deficiency-critical graph $G$ from $L_0=\{x_0\}$. If $x\in L_{N-1}$ satisfying $N(x)\cap L_N\neq \emptyset$, then there exists $y_1,y_2\in L_N$, such that $N(y_1)=N(y_2)=\{x\}$.	
\end{cor}

\begin{proof}
Applying  \cref{two level} to the last two levels $G[L_{N-1}\cup L_N]$ of $G$, i.e., $X=L_{N-1}$ and $Y=L_N\subseteq N(X)$, there is a matching $M$ of $G[L_{N-1}\cup L_N]$ satisfying:
\begin{enumerate}
	\item[(1)] $Y_M\subseteq N(X_M)$;
	\item[(2)] $Y_M$ is a stable set;
	\item[(3)] for any $x\in X_M$, there exist $y_1,y_2\in Y_M$, such that $N_{G_M}(y_1)=N_{G_M}(y_2)=\{x\}$,
\end{enumerate}
where $Y_M=L_N-V(M)$,
$X_M=\{x\in L_{N-1}-V(M):N(x)\cap Y_M\neq \emptyset\}$, and $G_M=G[L_{N-1}\cup L_N]-V(M)$. 
Let $G'=G-V(M)$. Each vertex of $G'$ can be connected to $x_0$ by the levelling as $Y_M\subseteq N(X_M)$. Thus $G'$ is still connected. Furthermore, $\kd(G')\ge \kd(G)$, since $G$ has a matching $M\cup M'$ missing exactly $\kd(G')$ vertices, where $M'$ is a maximum matching of $G'$. Since $G$ is deficiency-critical, we conclude that $G'=G$. 
Thus $M=\emptyset$,  $Y_M=L_N$, $X_M=\{x\in L_{N-1}:N(x)\cap L_N\neq \emptyset\}$, and $G_M=G[L_{N-1}\cup L_N]$. Hence for any vertex $x\in L_{N-1}$ satisfying $N(x)\cap L_N\neq \emptyset$, i.e., $x\in X_M$, according to (3), there exist $y_1,y_2\in L_N$, such that $N_{G_M}(y_1)=N_{G_M}(y_2)=\{x\}$. Therefore, $N(y_1)=N(y_2)=\{x\}$ as $N(y_i)\subseteq L_{N-1}\cup L_N$ for $i=1,2$.
\end{proof}

 When a vertex $x$ has two pendent edges, $xy$ and $xz$,   we refer to  this configuration as a {\em snail horn}. In this context, $x$ is called the {\em snail head} and $y, z$ are termed the corresponding {\em snail beards}. According to \cref{snail horn}, for nontrivial deficiency-critical graphs, we have the following corollary.
\begin{cor}\label{snail head}
Every nontrivial deficiency-critical graph $G$ contains a snail horn.
\end{cor}

The following algorithm constructs a large matching of $G$ and thus provides an upper bound of the deficiency of $G$.
\begin{algorithm}[H] 
\caption{Levelling-Matching (LM)} 
\label{alg:LM} 
\begin{algorithmic}[1] 
\REQUIRE A graph $G$ with levelling $V(G)=\bigcup\limits_{i=0}^N L_i$ from a snail head $x_0$. 
\ENSURE An upper bound of $\kd(G)$.
\STATE  Set $G_N=G$ and matching $M_{N+1}=M'_{N+1}=\emptyset$.
\WHILE{$1\leq i\le N$}

  \IF{$i\geq 1$}
      \STATE Applying \cref{two level} to the subgraph  $G_i[L_{i-1}\cup (L_i-V(M_{i+1}\cup M'_{i+1}))]$ induced by the last two levels of $G_i$, we can find a matching $M_i$ such that the sets 
\[Y_i:=L_i-V(M_{i+1}\cup M'_{i+1})-V(M_i) \text{ and } X_{i-1}:=\{u\in L_{i-1}-V(M_i):N(u)\cap Y_i\neq \emptyset\}\] satisfying 
\begin{enumerate}[(1)]
	\item $Y_i\subseteq N(X_{i-1})$,
	\item  $Y_i$ is a stable set,
	\item For any $u\in X_{i-1}$, there exist $u^+,u^-\in Y_i$, such that $N_G(u^+)\cap X_{i-1}=N_G(u^-)\cap X_{i-1}=\{u\}$.
	\item For any $v\in V(M_i)\cap L_{i-1}$, at most one vertex $u\in X_{i-1}$ does not satisfy that there exist $u^+,u^-\in Y_i$ such that $N_G(u^+)\cap (\{v\}\cup X_{i-1})=N_G(u^-)\cap (\{v\}\cup X_{i-1})=\{u\}$.
\end{enumerate}
Define 
\begin{align*}
M'_i:&=\{uu^+:u\in X_{i-1}\},\\
Z_i:&=Y_i-V(M'_i)=L_i-V(M_{i+1}\cup M'_{i+1}\cup M_i\cup M'_i),\\	
G_{i-1}:&=G_i-V(M_i\cup M'_i)-Z_i=\bigcup_{j=0}^{i-2}L_j\cup (L_{i-1}-V(M_i\cup M'_i)).
\end{align*}
\text{Update } $i:= i-1.$
  \ELSE
        \STATE Stop
    \ENDIF
\ENDWHILE
\RETURN 
$\kd(G)\le \Big|\bigcup_{i=1}^{N}Z_i\Big|=\sum_{i=1}^N |Z_i|$.

\end{algorithmic}
\end{algorithm}

\noindent{\bf Remark}: 
Note that \cref{two level} is valid for $G_i[L_{i-1}\cup (L_i-V(M_{i+1}\cup M'_{i+1}))]$ as $L_i-V(M_{i+1}\cup M'_{i+1})\subseteq L_i\subseteq N(L_{i-1})$. We claim that \cref{alg:LM} actually returns an upper bound of $\kd(G)$. 
Clearly, $L_i\subseteq Z_i\cup V(M_{i+1}\cup M'_{i+1}\cup M_i\cup M'_i)$ for any $i\in [1,N]$. Note that $L_0=\{x_0\}$ is saturated by $M'_1$ as $x_0$ is a snail head with two snail beards $x_0^+$ and $x_0^-$ in $L_1$, i.e., $x_0^+, x_0^-\in Y_1$.
Hence the matching $M:=\bigcup\limits_{i=1}^{N}(M_i\cup M'_i)$  can saturate all vertices of $G$ except $\bigcup\limits_{i=1}^{N}Z_i$.  Thus 
\[\kd(G)\le \Big|\bigcup_{i=1}^{N}Z_i\Big|=\sum_{i=1}^N |Z_i|.\]

\noindent{\bf Remark}: Throughout the following proofs, for any vertex $x\in X_{i-1}$, we consistently take $x^+,x^-\in Y_{i}$ satisfying $N(x^+)\cap X_i=N(x^-)\cap X_i=\{x\}$ given by conditions $(3)$ of \cref{alg:LM} unless stated otherwise.

\begin{lem}\label{Z_1}
$|Z_1|\le \alpha_l(G)-1$.	
\end{lem}
\begin{proof}
Note that $Y_1\cup \{x_0\}$ forms an induced star as $Y_1$ is a stable set. Thus $|Y_1|\le \alpha_l(G)$. Since $M'_1$ saturates a vertex of $Y_1$ and $x_0$, we obtain that  $|Z_1|=|Y_1|-1\le \alpha_l(G)-1$.
\end{proof}

\begin{lem}\label{|X_i-1|}
For $i>1$, $|Z_i|\le (\alpha_l(G)-2)|X_{i-1}|$.	
\end{lem}
\begin{proof}
Each $x\in X_{i-1}$ satisfies $|N(x)\cap Z_i|\le \alpha_l(G)-2$ since $\{x,x_{i-2},x^+\}\cup (N(x)\cap Z_i)$ is an induced star. Since $Z_i\subseteq Y_i\subseteq N(X_{i-1})$, we have 
\[|Z_i|=|Z_i\cap N(X_{i-1})|=\big|Z_i\cap \bigcup_{x\in X_{i-1}} N(x)\big| \le \sum_{x\in X_{i-1}} |N(x)\cap Z_i|\le (\alpha_l(G)-2)|X_{i-1}|.\]
\end{proof}

For $i\ge 1$, we say $L_{i}$ is {\it clean} if for any $u\in L_{i-1}$, $C(u)$ is a clique. $L_1$ is not clean as $x_0^+,x_0^-\in C(x_0)$ are non-adjacent.

\begin{lem}\label{messy level}
Let $i>1$. If $L_{i}$ is not clean, then $i\in A$.	
\end{lem}
\begin{proof}
Since $L_i$ is not clean, there exist non-adjacent vertices $v,w\in C(u)$ for some $u\in L_{i-1}$. Thus $V(uP_ux_0)\cup \{v,w,x_0^+,x_0^-\}$ forms bone $B_{i}$, which implies that $i\in A$. 
\end{proof} 

\begin{lem}\label{clean level no deficiency}
If $L_i$ is clean, then $|Z_i|=0$.	
\end{lem}
\begin{proof}
If $|Z_i|\neq 0$, then $|X_{i-1}|\neq 0$ by \cref{|X_i-1|}.
Thus there exist $u\in X_{i-1}$ and $u^+,u^-\in Y_i$ adjacent to $u$. $u^+u^-\not \in E(G)$ as $Y_i$ is a stable set.  This contradicts to the assumption that $L_i$ is clean. Hence $|Z_i|=0$.
\end{proof}

\begin{cor}\label{Z_inot=0}
If $|Z_i|\neq 0$, then $i\in A\cup \{1\}$.    
\end{cor}
\begin{proof}
Suppose to the contrary that $i\notin A\cup \{1\}$, thus $L_i$ is clean according to \cref{messy level}. Then $|Z_i|=0$ according to \cref{clean level no deficiency}, a contradiction.   
\end{proof}

The following lemma will be frequently used to find bones in a fundamental configuration. Note that if $P_u$ and $P_x$ are friendly at level $j$, then conditions $(1)$ and $(2)$ in \cref{friendly configuration} will be satisfied naturally.

\begin{lem}\label{friendly configuration}
Let $G$ be a deficiency-critical graph with snail head $x_0$ and admitting set $A$ consisting of odd integers. Let $V(G)=\bigcup_{i=0}^{N}L_i$ be the levelling of $G$ from $L_0=\{x_0\}$. 
Let $x\in L_{N-1}$ be another snail head.	
Let $v,w\in C(u)$ be non-adjacent for some $u\in L_{p-1}$ where $p\ge 3$.
Suppose $P_x$ and $P_u$ satisfies that 
\begin{enumerate}[(1)]
	\item $u_j=x_j$ for some $j\in [0,p-1]$,
	\item $u_kx_{k+1},x_ku_{k+1}\not \in E(G)$ for any $k>j$,
	\item There is no edge between $\{x_{p-1},x_p,x_{p+1}\}$ and $\{v,w\}$. 
\end{enumerate}
Then 
\begin{enumerate}[(i)]
	\item for any $k>j$, $u_kx_k\not \in E(G)$,
	\item $j+1\in A\cup \{1\}$ and $p-j,N-j,N+p-2j-1\in A$.
\end{enumerate}
 \end{lem}
 \begin{proof}
 (i) It can be directly checked that $V(x P_x x_0)\cup \{x^+,x^-,x_0^+,x_0^-\}$ forms a bone $B_N$, and $V(u P_u x_0)\cup \{v,w,x_0^+,x_0^-\}$ forms a bone $B_p$. Thus $N,p\in A$. Suppose, contrary to (i), that such a maximum integer $K=\max\{k>j:u_kx_k\in E(G)\}$ exists. Then $V(xP_xx_Ku_KP_uu)\cup \{v,w,x^+,x^-\}$ forms a bone $B_{N+p-2K}$. However, $N+p-2K$ is even as $N,p\in A$ are odd.
 This contradicts to that the admitting set $A$ contains only odd integers.

 (ii) If $j\neq 0$, then $V(x_jP_x x_0)\cup \{x_0^+,x_0^-,u_{j+1},x_{j+1}\}$ forms the bone $B_{j+1}$, and consequently, $j+1\in A$. If $j=0$, then $j+1=1$. Therefore, $j+1\in A\cup \{1\}$. 

Clearly, $N\ge p$ as $L_p\neq \emptyset$. Note that $j\neq p-1$, as otherwise, $x_{p-1}=u_{p-1}$ and it is adjacent to $v,w$, a contradiction. If $j\in [1,p-2]$, then $V(uP_uu_j)\cup \{x_{j-1},x_{j+1},v,w\}$ forms the bone $B_{p-j}$. Thus $p-j\in A$.  If $j=0$, then $p-j=p\in A$. In a word, $p-j\in A$.  
 If $j\neq 0$, then $V(xP_xx_j)\cup \{u_{j-1},u_{j+1},x^+,x^-\}$ forms the bone $B_{N-j}$, and consequently, $N-j\in A$. If $j=0$, then $N-j=N\in A$. Hence, $N-j\in A$. 
 Finally, as $V(xP_x x_j P_u u)\cup \{v,w,x^+,x^-\}$ forms the bone $B_{N+p-2j-1}$, we have $N+p-2j-1\in A$.
 \end{proof}

\section{Proof of \cref{main} }
\main*
\begin{proof}
For all connected induced subgraphs of $G$, we select $G'$ to have the maximum deficiency, and among those with the same deficiency, then the minimum order. Thus $\kd(G')\ge \kd(G)$.
\begin{claim}
$G'$ is deficiency-critical.	
\end{claim}
\begin{proof}
Given connected proper induced subgraph $G''$ of $G'$, we know that $G''$ is a connected induced subgraph of $G$ and $|V(G'')|<|V(G')|$. Thus $\kd(G'')<\kd(G')$ according to the definition of $G'$. Hence $G'$ is deficiency-critical.
\end{proof}
 
Hence we may assume $G$ is deficiency-critical as $G$ can be replaced by $G'$. If $G$ is trivial, then $\kd(G)\le |V(G)|\le 1$. Therefore, we may assume $G$ is nontrivial. According to \cref{snail head}, there is a snail horn in $G$, i.e., vertices $x_0,x_0^+,x_0^-$ satisfying $N(x_0^+)=N(x_0^-)=\{x_0\}$.
Let $V(G)=\bigcup\limits_{i=0}^N L_i$ be the levelling of $G$ from $x_0$. Applying \cref{alg:LM} to $G$ with levelling   $V(G)=\bigcup\limits_{i=0}^N L_i$, we have 
$X_{i-1}$, $Y_i$ and $Z_i$ for $1\le i\le N$, and 
\[\kd(G)\le \Big|\bigcup_{i=1}^{N}Z_i\Big|=\sum_{i=1}^N |Z_i|.\]

Now all we need is to bound $\sum\limits_{i=1}^N |Z_i|$. By \cref{Z_1}, $|Z_1|\le \alpha_l(G)-1\le n-2$. For any even $i\ge 2$, $|Z_i|=0$ by \cref{Z_inot=0} as $i\notin A\cup \{1\}$. When we want to bound $|Z_i|$ generally, we fix $i$ and regard $G_i$ as the host graph. For convenience, let $L'_i=Y_i$ and $L'_{i-1}=X_{i-1}$. For $j=i-2,\cdots ,1$, recursively let $L'_j$ be a minimal subset of $L_j$ satisfying that $L'_{j+1}\subseteq N(L'_j)$. By the minimality of $L'_j$, for any $u\in L'_j$, there exists at least one private neighbor $p(u)\in L'_{j+1}$ satisfying $N(p(u))\cap L'_j=\{u\}$. Hence we can define a private neighbor function 
\[p:\bigcup_{j=1}^{i-2}L'_j\longrightarrow \bigcup_{j=2}^{i-1}L'_j\] 
by fixing a private neighbor of $u$ for any $u\in \bigcup\limits_{j=1}^{i-2}L'_j$. Let $p^{(0)}(u)=u$ and iteratively define $p^{(j)}(u)=p(p^{(j-1)}(u))$. For any distinct $u,v\in L'_j$, $p^{(k-1)}(u)p^{(k)}(v), \;p^{(k-1)}(v)p^{(k)}(u)\not\in E(G)$ according to the definition of the private neighbor function $p$.
Note that $L'_j$ is temporary for $j\in [i]$, as they depend on the host graph $G_i$.

\begin{claim}\label{L' stable}
For fixed $i\in[3,N]$, $L'_j$ is a stable set for  $j\in [i]$.	
\end{claim}
\begin{proof}
Suppose not and let $j$	be the maximum such that $L'_j$ is not stable. Then $j\neq i$ as $L'_i=Y_i$ is stable. Thus there exist vertices $u,v\in L'_j$ satisfying $uv\in E(G)$. Furthermore, for any $k\in [1,i-j-1],\; p^{(k)}(u)p^{(k)}(v)\not \in E(G)$ as $L'_{j+k}$ is stable. Thus 
\begin{align*}
&\{p^{(0)}(u),\cdots , p^{(i-j-1)}(u), p^{(i-j-1)}(u)^+, p^{(i-j-1)}(u)^-\}\\
&\cup\{p^{(0)}(v),\cdots ,p^{(i-j-1)}(v), p^{i-j-1)}(v)^+, p^{(i-j-1)}(v)^-\}	
\end{align*}
forms the even bone $B_{2i-2j}$, which is not admitted in $G_i$, a contradiction. 
\end{proof}

\begin{claim}\label{L'_1}
$|L'_1|\le n-3$.	
\end{claim}
\begin{proof}
Each vertex $u\in L'_1$ has a private neighbor $p(u)\in L'_2\subseteq L_2$. Thus $x_0^+,x_0^-\not\in L'_1$ as $N(x_0^+)=N(x_0^-)=\{x_0\}$. Therefore, $\{x_0,x_0^+,x_0^-\}\cup L'_1$ forms an induced star, implying that $|L'_1|\le n-3$.	
\end{proof}

\begin{claim}\label{ditui} 
$|L'_{j+1}|\le (n-2)|L'_{j}|$ for $j\in [i-1]$.	
\end{claim}
\begin{proof}
For any $u\in L'_j$, choose a neighbor $u_{j-1}$ of $u$ in $L_{j-1}$, then $\{u,u_{j-1}\}\cup (N(u)\cap L'_{j+1})$ forms an induced star as $L'_{j+1}$ is a stable set by \cref{L' stable}. 
Hence $|N(u)\cap L'_{j+1}|\le n-2$ as $G$ is $K_{1,n}$-free.
Since $L'_{j+1}\subseteq N(L'_j)=\bigcup_{u\in L'_j} N(u)$, 
\[|L'_{j+1}|=\Big|\bigcup_{u\in L'_j}(N(u)\cap L'_{j+1})\Big| \le \sum_{u\in L'_j}|N(u)\cap L'_{j+1})|\le (n-2)|L'_j|.\]
\end{proof}

\begin{claim}\label{pintui}
If $L_{j+1}$ is clean, then $|L'_{j+1}|=|L'_{j}|$.  	
\end{claim}
\begin{proof}
For any $u\in L'_j$, $N(u)\cap L'_{j+1}$ is a clique as $L_{j+1}$ is clean, also a stable set as $L'_{j+1}$ is stable. Hence $N(u)\cap L'_{j+1}=\{p(u)\}$ is of size $1$. Therefore, $L'_{j+1}=\{p(u):u\in L'_j\}$. Thus $|L'_j|=|L'_{j+1}|$.	
\end{proof}

\begin{claim}
For any odd $i>1$, $|Z_i| \le (n-3)^2 (n-2)^{\frac{i-3}{2}}$.  	
\end{claim} 

\begin{proof}
\begin{align*}
|Z_i|
&=|Y_i|-|X_{i-1}|=|L'_i|-|L'_{i-1}|	&&\text{ as } |V(M'_i)\cap Y_i|=|X_{i-1}|\\
&\le (n-3)|L'_{i-1}| &&\text{ by }\cref{ditui}\\
&=(n-3)|L'_{i-2}| &&\text{ by } \cref{messy level} \text{ and } \cref{pintui}\\
&\le (n-3)(n-2)|L'_{i-3}| = (n-3)(n-2)|L'_{i-4}|\\ 
&\;\vdots\\
&\le (n-3)(n-2)^{\frac{i-3}{2}}|L'_2|=(n-3)(n-2)^{\frac{i-3}{2}}|L'_1|\\
&\le (n-3)^2 (n-2)^{\frac{i-3}{2}}&&\text{ by }\cref{L'_1}.
\end{align*}
\end{proof} 

\begin{claim}
If $i\ge m$, then $|Z_i|\le (n-3)(n-2)^{\frac{m-3}{2}}$.	
\end{claim}
\begin{proof}
We may assume $i$ is odd since $|Z_i|=0$ for any even $i$.
We first prove that $|L'_{i-m+1}|=1$. Suppose not, then there exist distinct $u,v\in L'_{i-m+1}$. Let $P_u$ and $P_v$ be friendly at level $j$, where $j\in [0,i-m]$. Thus for any $k\in [j+1,i-m]$, $u_kv_{k+1},u_{k+1}v_k\not \in E(G)$. 
Let $P_{p^{(i-j-1)}(u)}=p^{(i-j-1)}(u)p^{(i-j-2)}(u)\cdots p^1(u)uP_ux_0$ and $P_{p^{(i-j-1)}(v)}=p^{(i-j-1)}(v)p^{(i-j-2)}(v)\cdots p^1(v)vP_{v}x_0$. In graph $G_i$, $p^{(i-j-1)}(u)$ is a snail head with beards $p^{(i-j-1)}(u)^+,p^{(i-j-1)}(u)^-$, and 
\[p^{(i-j-1)}(v)^+,p^{(i-j-1)}(v)^-\in N(p^{(i-j-1)}(v))\cap L_i\]
are non-adjacent. 
Applying \cref{friendly configuration} to  $P_{p^{(i-j-1)}(u)}$ and  $P_{p^{(i-j-1)}(v)}$, we have $i-j,2i-2j-1\in A$. However, $i-j\ge m$ as $j\le i-m$.
This leads to a contradiction as $2i-2j-1=(i-j)+(i-j)-1\in A$. Thus $|L'_{i-m+1}|=1$. Therefore,
\begin{align*}
|Z_i|
&=|Y_i|-|X_{i-1}| =|L'_i|-|L'_{i-1}|\\
&\le (n-3)|L'_{i-1}| =  (n-3)|L'_{i-2}| \\ 
&\;\vdots\\
&\le (n-3)(n-2)^{\frac{m-3}{2}}|L'_{i-m+1}|,\\
&\le (n-3) (n-2)^{\frac{m-3}{2}}.
\end{align*}	
\end{proof}

If the level number $N< 3m-4$, then
\begin{align*}
\sum_{i=1}^{N}|Z_i|&=|Z_1|+ \sum_{\text{odd }i\in [3,m-2]}|Z_i|+\sum_{\text{odd }i\in [m,N]}|Z_i|\\
&\le (n-2)+\sum_{i=0}^{\frac{m-5}{2}} (n-3)^2(n-2)^i+(m-2)(n-3)(n-2)^{\frac{m-3}{2}}\\
&=(m-1)(n-3)(n-2)^{ \frac{m-3}{2}}+1.	
\end{align*}

Now we assume that $N\ge 3m-4$. According to \cref{two level}, there is a snail head $x$ in $L_{N-1}$ with beards $x^+,x^-$. Note that $V(xP_xx_0)\cup \{x^+,x^-,x_0^+,x_0^-\}$ forms the bone $B_N$. Thus $N\in A$. 

\begin{claim}\label{middle level clean}
For any $p\in [2m-4,N-m]$, $|Z_p|=0$.	
\end{claim}
\begin{proof}
Suppose $|Z_p|\neq 0$ for some $p\in [2m-4,N-m]$, then $X_{p-1}\neq \emptyset$ by \cref{|X_i-1|}. Thus there exists $u\in X_{p-1}$ and $u^+,u^-\in Y_p\cap N(x)$ such that $u^+u^-\notin E(G)$. 
By Corollary~\ref{Z_inot=0}, $p\in A\cup \{1\}$.  
Thus $p\in [2m-3,N-m-1]$.
Let $P_u$ and $P_x$ be friendly at level $j$, where $j\in [0,p-1]$. 
We check that there is no edge between $\{x_{p-1},x_p,x_{p+1}\}$ and $\{u^+,u^-\}$ first. 

For $i\in \{p-1,p,p+1\}$, let $a_i=|N(x_{i})\cap \{u^+,u^-\}|$. We show that $a_i=0$ for $i=p+1,p,p-1$ one by one. If $a_i=1$ for some $i\in \{p-1,p,p+1\}$, we may assume $N(x_{i})\cap \{u^+,u^-\}=\{u^+\}$ for convenience in the following.

If $a_{p+1}=2$, then $V(x P_x x_{p+1})\cup \{u^+,u^-,x^+,x^-\}$ forms the bone $B_{N-p-1}$.
This leads to a contradiction as $N=(N-p-1)+p+1$ while $p,N-p-1\in A_{\ge m}$. If $a_{p+1}=1$,  then $V(x P_x x_{p+1} u^+u)\cup \{u^-,u_{p-2},x^+,x^-\}$ forms the bone $B_{N-p+1}$, a contradiction again as $N=(N-p+1)+p-1$ while $N-p+1, p\in A_{\ge m}$. Thus $a_{p+1}=0$. 

If $a_{p}=2$, then $V(x P_x x_{p})\cup \{u^+,u^-,x^+,x^-\}$ forms an even bone, a contradiction. Now assume $a_p=1$. If $x_pu\not \in E(G)$, then $V(x P_x x_p u^+ u)\cup \{u_{p-2},u^-,x^+,x^-\}$ forms an even bone, a contradiction. On the other hand, if $x_pu\in E(G)$, then $V(x P_x x_p u)\cup \{u_{p-2},u^-,x^+,x^-\}$ forms the bone $B_{N-p+1}$, a contradiction. Thus $a_p=0$.

If $a_{p-1}\not=0$, then $x_{p-1}u^+\in E(G)$. Then $V(xP_xx_{p-1})\cup \{x^+,x^-, u^+,x_{p-2}\}$ forms the bone $B_{N-p+1}$, a contradiction too. Thus $a_{p-1}=0$. Hence there is no edge between $\{x_{p-1},x_p,x_{p+1}\}$ and $\{u^+,u^-\}$.

According to \cref{friendly configuration}, $j+1,N+p-2j-1, p-j,N-j\in A\cup \{1\}$.
If $j\ge m-1$, then $j+1\in A_{\ge m}$. $N-j\in A_{\ge m}$ as $N-j\ge N-(p-1)> m>1$.
Therefore, $N=(j+1)+(N-j)-1 \in A_{\ge m}$, a contradiction.

On the other hand, if $j<m-1$. Then $j\le m-3$ as $j$ is even and $m$ is odd.
Note that $N-j>p-j\ge (2m-3)-(m-3)=m$. Thus $p-j,\,N-j\in A_{\ge m}$. This leads to a contradiction as $N+p-1-2j=(p-j)+(N-j)-1$.
The proof of \cref{middle level clean} is complete.
\end{proof}

In summary, we have the upper bound of $\sum_{i=1}^{N}|Z_i|$   in this case as follows:
\begin{align*}
\sum_{i=1}^{N}|Z_i|&=|Z_1|+ \sum_{\text{odd }i\in [3,m-2]}|Z_i|+\sum_{\text{odd }i\in [m,2m-5]\cup [N-m+1,N]}|Z_i|\\
&\le (n-2)+\sum_{i=0}^{\frac{m-5}{2}} (n-3)^2(n-2)^i+ (m-1)(n-3)(n-2)^{\frac{m-3}{2}}\\
&=m(n-3)(n-2)^{ \frac{m-3}{2} }+1.	
\end{align*}
Therefore, $\kd(G)\le \sum_{i=1}^{N}|Z_i|\le m(n-3)(n-2)^{ \frac{m-3}{2} }+1$.

\bigskip
When $m=3$, we have $A_{\ge m}=A$ as $2\not \in A$.
\begin{claim}\label{Y_N}
$Y_N\subseteq N(x)$.	
\end{claim} 
\begin{proof}
Suppose $Y_N\not \subseteq N(x)$, then $X_{N-1}\neq \{x\}$ as $Y_N\subseteq N(X_{N-1})$. Thus there exists a vertex $u\in X_{N-1}$, and $u^+,u^-\in Y_N$, such that $N(u^+)=N(u^-)=\{u\}$. Let $P_x$ and $P_u$ be friendly at level $j$, where $j\in [0,N-2]$. 
There is no edge between $x=x_{N-1},x_N,x_{N+1}$ and $u^+,u^-$ as $N(u^+)=N(u^-)=\{u\}$. 
Applying \cref{friendly configuration} to $P_u$ and $P_x$, we have $N+p-2j-1, N-j\in A$ where $p=N$. This leads to a contradiction as $2N-2j-1=(N-j)+(N-j)-1$. Thus $Y_N\subseteq N(x)$.	
\end{proof}
If $N<3m-4=5$, then $N=3$ as $N>1$ is odd. $|Z_3|\le |Y_3|-1\le n-3$ as $\{x,x_1\}\cup Y_3$ is an induced star.
$\kd(G)=|Z_1|+|Z_3|\le 2n-5$.
On the other hand, if $N\ge 5$, by \cref{middle level clean}, $|Z_i|=0$ for $i\in [2,N-3]$. 

\begin{claim}\label{Y_N-2}
If $N\ge 5$, then $Y_{N-2}\subseteq N(x)$. 
\end{claim}
\begin{proof}
Suppose $Y_{N-2}\not \subseteq N(x)$, then there exists $w\in Y_{N-2}$ nonadjacent to $x$. Since $Y_{N-2}\subseteq X_{N-3}$, there exists $u\in X_{N-3}$ adjacent to $w$. Since $|N(u)\cap Y_{N-2}|\ge 2$, there exists $v\in Y_{N-2}\cap N(u)$ distinct from $w$. As $Y_{N-2}$ is a stable set, $vw\not \in E(G)$. Since $V(uP_ux_0)\cup \{x_0^+,x_0^-,v,w\}$ forms the bone $B_{N-2}$, we have $N-2\in A$. 
Let $P_x$ and $P_u$ be friendly at level $j$, where $j\in [0,N-3]$. 

We check that there is no edge between $\{x_{N-3},x_{N-2}, x=x_{N-1}\}$ and $\{v,w\}$.
For $i=N-3, N-2, N-1$, let $a_i=|N(x_i)\cap \{v,w\}|$.
If $xv\in E(G)$, then $\{x,v,u,u_{N-4},w,x^+,x^-\}$ forms the bone $B_3$, a contradiction as $N=(N-2)+3-1\in A$. Thus $a_{N-1}=0$. If $a_{N-2}=2$, then $\{x,x_{N-2},v,w,x^+,x^-\}$ forms the bone $B_2$, a contradiction. Now suppose $a_{N-2}=1$, i.e. assume $x_{N-2}v\in E(G)$.
If $x_{N-2}u\not \in E(G)$, then $\{x, x_{N-2},v, u, u_{N-4},w,x^+,x^-\}$ forms the bone $B_4$, a contradiction. If $x_pu\in E(G)$, then $\{x, x_{N-2}, u, u_{N-4},w,x^+,x^-\}$ forms the bone $B_3$, a contradiction again. Thus $a_{N-2}=0$.
If $a_{N-3}\not=0$, assume $x_{N-3}v\in E(G)$. Then $\{x,x_{N-2}, x_{N-3},x^+,x^-, v,x_{N-4}\}$ forms bone $B_3$, a contradiction too. Thus $a_{N-3}=0$. Hence there is no edge between $\{x_{N-3},x_{N-2},x_{N-1}\}$ and $\{v,w\}$. By \cref{friendly configuration}, $N+p-2j-1, p-j,N-j\in A$ where $p=N-2$ and $j=\flev(u,x)$. This leads to a contradiction as $2N-2j-3=(N-j)+(N-2-j)-1$. Thus $Y_N\subseteq N(x)$. 
\end{proof}
Hence $\{x\}\cup Y_{N-2}\cup Y_N$ is an induced star by Claims \ref{Y_N} and \ref{Y_N-2}, implying that $|Y_{N-2}\cup Y_N|\le n-1$. 
Therefore 
\begin{align*}
\kd(G)&=|Z_1|+|Z_{N-2}|+|Z_N|\\
&\le (n-2)+(|Y_{N-2}|-1)+(|Y_N|-1) && \text{ as } X_{N-3},X_{N-1}\neq \emptyset \\
&=(n-2)+|Y_{N-2}\cup Y_N|-2\\ 
&\le 2n-5
\end{align*}
The proof of \cref{main} is complete.
\end{proof}


\section{Proofs of Theorems \ref{two odd bone} and \ref{only one even bone}}
During this section, we may always assume that $G$ is a deficiency-critical graph with snail horn $x_0,x_0^+,x_0^-$. Applying \cref{alg:LM} to the levelling $V(G)=\bigcup_{i=0}^N L_i$ of $G$ from $x_0$, we can get the corresponding sets $X_{i-1},Y_i$ and $Z_i$ for $i\in [N]$. By Lemma~\ref{Z_1}, $|Z_1|\le \alpha_l(G)-1\le n-2$.
If $N=1$, then $\kd(G)= |Z_1|\le \alpha_l(G)-1\le n-2$. Therefore, we may always assume $N>1$. Thus there exists $x\in X_{N-1}$ with beards $x^+,x^-$. Let $X_{N-1}=\{x^i:i\in [s]\}$, where $x=x^1, s=|X_{N-1}|$. Let $A$ be the admitting set of $G$. Then $N\in A$ since $V(xP_xx_0)\cup \{x_0^+,x_0^-,x^+,x^-\}$ forms the bone $B_N$.  By \cref{Z_inot=0}, $|Z_i|=0$ for any $i\not \in \{1\}\cup A$ .
 
\twooddbone*

\begin{proof}
Recall that $N\in A=\{p,q\}$, where $q=2p\pm 1>p$. 

Assume $N=p$ first. 
If $s=1$, then by \cref{|X_i-1|}, $\kd(G)\le |Z_1|+|Z_p|\le (n-2)+(n-3)<3n-8< n^2-3n+1$. This is exactly what we need to prove. 
Now assume $s>1$. Let $P_x$ and $P_{x^2}$ be friendly at level $j$, where $j\in [0,p-2]$. There is no edge between $x^2_{p-1}$ and $\{x^+,x^-\}$. By \cref{friendly configuration},  $j+1,2p-2j-1\in \{1,p,q\}$. Thus $j=0$ and consequently, $2p-2j-1=2p-1\in \{1,p,q\}$ and $q=2p-1$. By symmetry, $\flev(x^i,x^k)=0$ for any different $i,k\in [s]$.  
Thus $P_{x^i}$ and $P_{x^k}$ are automatically friendly at level $0$ for any distinct $i,k\in [s]$. 
According to \cref{friendly configuration}, $x^i_1x^k_1\not \in E(G)$. Therefore, $\{x_0,x_0^+,x_0^-,x^i_1:i\in [s]\}$ forms an induced star $K_{1, s+2}$. Consequently, $s\le n-3$.
Hence $\kd(G)\le (n-2)+(n-3)s< n^2-3n+1$. 

Now assume $N=q$. Thus $x\in X_{N-1}=X_{q-1}$. Consider first the case that $q=2p+1$. If $s>1$, then by \cref{friendly configuration}, $j+1,2q-2j-1,q-j\in \{1,p,q\}$, where $j=\flev(x,x^2)\in [0,q-2]$. $2q-2j-1\le q$ implies that $j\ge p$. Then $j+1\not\in \{1,p,q\}$, a contradiction. Hence $s=1$ and $|Z_q|\le n-3$ by \cref{|X_i-1|}. If $|X_{p-1}|>1$, let $u,v\in X_{p-1}$ be distinct.
Applying \cref{friendly configuration} in graph $G_p$ (Just to ensure that $u$ is a snail head of $G_p$), $j+1,2p-2j-1,p-j\in \{1,p,q\}$, where $j=\flev(u,v)\in [0,p-2]$. $j+1\in \{1,p,q\}$ implies that $j=0$, then $2p-2j-1=2p-1\not\in \{1,p,q\}$, a contradiction. Hence $|X_{p-1}|\le 1$ and thus $|Z_p|\le (\alpha_l(G)-2)|X_{p-1}|\le n-3$. Therefore, $\kd(G)=|Z_1|+|Z_p|+|Z_q|\le 3n-8$.

Now we consider the case that $q=2p-1$.
\begin{claim}
$|Z_{p}|\le n-3$.	
\end{claim}
\begin{proof}
Suppose not, then $|X_{p-1}|>1$ according to \cref{|X_i-1|}. Thus there exist distinct $u,v\in X_{p-1}$  as snail heads in $G_p$.  Applying \cref{friendly configuration} to $G_p$, we have $\flev(u,v)+1\in \{1,p,q\}$ (thus $\flev(u,v)=0$ as $\flev(u,v)<p-1$) and $uP_ux_0P_v v$ is an induced path. Any selected $P_u$ and $P_v$ for $u,v$ are friendly at level $0$. 
Hence we may assume $P_x$ and $P_u$ are friendly in the meantime. According to $(4)$ of \cref{alg:LM}, we know $x_{p-1}$ can not destroy both $u$ and $v$ as snail heads of $G_p$. Thus we may assume $N(x_{p-1})\cap \{u^+,u^-\}=\emptyset$,  where $u^+,u^-$ are beards of $u$, and $v^+,v^-$ are beards of $v$ in $G_p$.
For $i\in \{p+1,p\}$, let $b_i=|N(x_{i})\cap \{u^+,u^-,v^+,v^-\}|$. We will prove that $b_i=0$ for $i=p+1,p$ one by one. $x_{p-1}$ is not necessary for this proof. Therefore, if $b_i=1$ for some $i\in \{p,p+1\}$, we may assume $N(x_{i})\cap \{u^+,u^-,v^+,v^-\}=\{u^+\}$ for convenience in the following.

If $b_{p+1}\ge 2$, then the union 
\[V(x P_x x_{p+1})\cup \{x^+,x^-\} \cup \big(N(x_{p+1})\cap \{u^+,u^-,v^+,v^-\}\big)\]
necessarily contains an even bone $B_{p-2}$, a contradiction. If $b_{p+1}=1$, then $$V(x P_x x_{p+1}u^+u P_u x_0 P_v v)\cup \{x^+,x^-,v^+,v^-\}$$ forms the bone $B_{q+p-1}$, a contradiction. Therefore, $b_{p+1}=0$.

If $b_{p}\ge 2$, then the union 
\[V(x P_x x_p)\cup \{x^+,x^-\} \cup \big(N(x_{p})\cap \{u^+,u^-,v^+,v^-\}\big)\]
contains an even bone $B_{p-2}$, yielding a contradiction to the admitting set $A=\{p, q\}$. If $b_p=1$, then $x_p v\not \in E(G)$; otherwise, $V(x P_x x_p)\cup \{x^+,x^-,u^+,v\}$ forms the bone $B_{p-1}$, a contradiction.
Hence the union $V(x P_x x_p  u P_{u} x_0 P_{v}v)\cup \{x^+,x^-,v^+,v^-\}$ 
or $V(x P_x x_p u^+ u P_{u} x_0 P_{v}v)\cup \{x^+,x^-,v^+,v^-\}$ forms a bone larger than $B_q$, according to whether $x_pu\in E(G)$, a contradiction. Therefore, $b_p=0$.

Then there is no edge between $\{x_{p-1},x_p,x_{p+1}\}$ and $\{u^+,u^-\}$. 
By \cref{friendly configuration}, $j+1,p+q-2j-1\in \{1,p,q\}$ where $j=\flev(x,u)<p-1$ as $u\neq x_{p-1}$. Thus $j=0$, this leads to a contradiction as $p+q-2j-1=p+q-1\not \in \{1,p,q\}$.
\end{proof}
\begin{claim}
$|Z_q|\le (n-3)(n-2)$.	
\end{claim}
\begin{proof}
Suppose not, then $s=|X_{q-1}|>n-2$ according to \cref{|X_i-1|}. Let $P_x$ and $P_{x^2}$ be friendly at level $\flev(x,x^2)$. By \cref{friendly configuration}, $\flev(x,x^2)+1,q+p-2\flev(x,x^2)-1\in \{1,p,q\}$. Thus $\flev(x,x^2)=p-1$. By symmetry, $\flev(x^i,x^k)=p-1$ for any distinct $i,k\in [s]$. 

We will show that there exist $P_{x^1},P_{x^2},\cdots ,P_{x^i}$ for any $i\in [2,s]$
such that $x^{\ell}_{p-1}$ are all the same for each $\ell\in [1,i]$ by induction on $i$. For $i=2$, $P_x$ and $P_{x^2}$ fulfill requirement. We now handle the $i$ case based on the $i-1$ case. Thus there exist $P_{x^1},P_{x^2},\cdots ,P_{x^{i-1}}$ such that $x^{\ell}_{p-1}$ are all the same for each $\ell\in [1,i-1]$. 
Let $Q_{x^i}=x^i x^i_{q-2}\cdots x^i_px^i_{p-1}$ be a shortest path connecting $x^i$ and $L_{p-1}$. Thus $x^i_p\in L_p$ and $x^i_{p-1}\in L_{p-1}$. If $x^i_p x_{p-1}\in E(G)$, then let $P_{x^i}=x^iQ_{x^i}x^i_px_{p-1}P_x x_0$, thus $x^i_{p-1}=x_{p-1}$. $P_{x^1},P_{x^2},\cdots ,P_{x^{i-1}}$ and $P_{x^i}$ fulfill requirement. Otherwise, we may assume that $x^i_p x_{p-1}\notin E(G)$ now.
we claim that $x^i_{p-1}x^{\ell}_{p}\in E(G)$ for any $\ell\in [1,i-1]$. 
Suppose to the contrary that $x^i_{p-1}x^{\ell}_{p}\not \in E(G)$ for some $\ell$. 
Thus $x^i_k x^{\ell}_{k+1},x^i_{k+1}x^{\ell}_k\not \in E(G)$ for any $k\ge p$ as $\flev(x^i,x^\ell)=p-1$. Let $P_{x^i_{p-1}}$ and $P_{x^{\ell}_{p-1}}$ be friendly at level $j$, where $j=\flev(x^{\ell}_{p-1},x^i_{p-1})<p-1$. Let $P_{x^i}=x^iQ_{x^i}x^i_{p-1} P_{x^i_{p-1}}x_0$ and $P'_{x^{\ell}}=x^{\ell}P_{x^{\ell}} x^{\ell}_{p-1} P_{x^{\ell}_{p-1}}x_0$. Applying \cref{friendly configuration} to $P_{x^i}$ and $P'_{x^{\ell}}$, we have $2q-2j-1\in \{p,q\}$. Thus $2q-2j-1\le q$, which contradicts to $j<p-1$. Hence $x^i_{p-1}x^{\ell}_{p}\in E(G)$ for any $\ell\in [1,i-1]$. Update $P_{x^{\ell}}$ as $x^{\ell}P_{x^{\ell}}x^{\ell}_px^i_{p-1}P_{x^i_{p-1}}x_0$ for all $\ell\in [1,i-1]$, then $x^{\ell}_{p-1}=x^i_{p-1}$ for all $\ell\in [1,i-1]$. Thus $P_{x^i}$ and the updated $P_{x^1},P_{x^2},\cdots ,P_{x^{i-1}}$ fulfill requirement. This completes the induction. Hence there exist $P_{x^1},P_{x^2},\cdots ,P_{x^s}$ such that $x^{\ell}_{p-1}$ are all the same for each $\ell\in [1,s]$. 

For any $i,k\in [s]$, $P_{x^i}$ and $P_{x^k}$ are already friendly at level $p-1$ as $x^i_{p-1}=x^k_{p-1}$ and $\flev(x^i,x^k)=p-1$.
Applying \cref{friendly configuration} to $P_{x^i}$ and $P_{x^k}$, we have $x^i_p x^k_p\not \in E(G)$. Thus $\{x_{p-1},x_{p-2},x^i_p:i\in [1,s]\}$ forms an induced star, which implies $\alpha_l(G)\ge s+1>n-1$, a contradiction.
\end{proof}
Therefore, $\kd(G)=|Z_1|+|Z_p|+|Z_q|\le (n-2)+(n-3)+(n-3)(n-2)=n^2-3n+1$. 
Now we characterize the deficiency-critical graphs attaining the extremum, respectively. Either $q=2p-1$ or $q=2p+1$, the deficiency-critical graphs $G$ satisfy that $|Z_1|=n-2$ and $|Y_1|=n-1$. We have the following claim in both cases.
\begin{claim}
There exists a stable set $S\subseteq L_1$ satisfying $|S|=n-1$ and $|N(x_2)\cap S|=1$.     
\end{claim}
\begin{proof}
If $|N(x_2)\cap Y_1|\ge 2$, then $V(xP_xx_2)\cup \{x^+,x^-\}\cup (N(x_2)\cap Y_1)$ contains a bone $B_{q-2}$, a contradiction, unless $p=3,q=5$. In this case, $\{x_2,x^i_3:i\in [n-2]\}\cup (N(x_2)\cap Y_1)$ is an induced star, contradicting to the assumption that $G$ is $K_{1,n}$-free. Therefore, $|N(x_2)\cap Y_1|<2$.
If $|N(x_2)\cap Y_1|=1$, then we take $S=Y_1$. Thus we may assume $N(x_2)\cap Y_1=\emptyset$. Hence $x_1\notin Y_1$. 

If $|N(x_1)\cap Y_1|\ge 2$, then  $V(xP_xx_2x_1)\cup \{x^+,x^-\}\cup (N(x_1)\cap Y_1)$ contains a bone $B_{q-1}$, a contradiction. If $|N(x_1)\cap Y_1|=0$, then $\{x_0,x_1\}\cup Y_1$ forms a star $K_{1,n}$, a contradiction too. Hence $|N(x_1)\cap Y_1|=1$. Let $S=Y_1\cup \{x_1\}\setminus (N(x_1)\cap Y_1)$. Thus $S$ is a stable set of size $n-1$ and $N(x_2)\cap S=\{x_1\}$, as desired.    
\end{proof}

One can check directly that $S_{n-1}^p$ is a $K_{1,n}$-free deficiency-critical graph with admitting set $\{p,2p-1\}$ satisfying $\kd(G)=n^2-3n+1$ from the construction. We show that $S_{n-1}^p$ is the only deficiency-critical graph attaining this extremum. Let $G$ be a $K_{1,n}$-free deficiency-critical graph with admitting set $\{p,2p-1\}$ satisfying $\kd(G)=n^2-3n+1$. Tracing the condition to realize $\kd(G)=n^2-3n+1$ in the above proof, we have $N=2p-1$, $|Z_{2p-1}|=(n-3)(n-2)$ and thus $|X_{2p-2}|=n-2$. For all $x^s\in X_{2p-2}$, $x^s_{p-1}=x_{p-1}$ is the same vertex. Therefore, $G$ contains $S_{n-1}^p$ as an induced subgraph. Since $G$ is deficiency-critical and $\kd(G)=\kd(S_{n-2}^p)$, graph $G$ must be $S_{n-2}^p$. Therefore, $S_{n-2}^p$ is the unique deficiency-critical graph attaining this extremum. 

Similarly, $T_{n-2}^p$ is a $K_{1,n}$-free deficiency-critical graph with admitting set $\{p,2p+1\}$ satisfying $\kd(T_{n-2}^p)=3n-8$. Let $G$ be a $K_{1,n}$-free deficiency-critical graph with admitting set $\{p,2p+1\}$. Considering the condition to ensure that $\kd(G)=3n-8$, we have $|Z_p|=|Z_{2p+1}|=n-3$, $|X_{p-1}|=|X_{2p}|=1$. Thus $|Y_p|=|Y_{2p+1}|=n-2$. Let $X_{p-1}=\{u\}$ and $Q_x=xx_{2p}\cdots x_{p+1}$ be a shortest path connecting $x$ and $L_{p+1}$. 
\begin{claim}
$Y_p\subseteq N(x_{p+1})$.    
\end{claim}  
\begin{proof}
Assume there exists $x_p\in Y_p\cap N(x_{p+1})$ first. Thus if there exists any $y\in Y_p\setminus N(x_{p+1})$, $V(xQ_xx_{p+1}x_pu)\cup \{x^+,x^-,y, u_{p-2}\}$ forms bone $B_{p+2}$, a contradiction. Hence $Y_p\subseteq N(x_{p+1})$ at this time. 

Now suppose to the contrary that $Y_p\cap N(x_{p+1})=\emptyset$. Take $x_p\in L_p$ adjacent to $x_{p+1}$, thus $x_p\notin Y_p$. If $|N(x_p)\cap Y_p|\ge 2$, then $V(xQ_xx_{p+1}x_p)\cup \{x^+,x^-,\}\cup (N(x_p)\cap Y_p)$ contains bone $B_{p+1}$, a contradiction. If $x_py\in E(G)$ for exactly one $y\in Y_p$, then $V(xQ_xx_{p+1}x_pyu)\cup \{x^+,x^-,u_{p-2},y'\}$ forms bone $B_{p+3}$, where $y'\in Y_p\setminus \{y\}$ or $V(xQ_xx_{p+1}x_pu)\cup \{x^+,x^-,u_{p-2},y'\}$ forms bone $B_{p+2}$ according to whether $x_pu\in E(G)$. This leads to a contradiction as $p+2,p+3\notin A$. Hence $N(x_p)\cap Y_p=\emptyset$. Take $x_{p-1}\in L_{p-1}$ adjacent to $x_p$.
$x_pu\notin E(G)$ otherwise $V(xQ_xx_{p+1}x_pu)\cup \{x^+,x^-,u^+,u^-\}$ forms bone $B_{p+2}$. Thus $x_{p-1}\neq u$. $N(x_{p-1})\cap Y_p=\emptyset$ otherwise $V(xQ_xx_{p+1}x_px_{p-1})\cup \{x^+,x^-,x_{p-2},y\}$ forms bone $B_{p+2}$, where $y\in N(x_{p-1})\cap Y_p$, and $x_{p-2}\in N(x_{p-1})\cap L_{p-2}$. Let $P_{x_{p-1}}$ and $P_u$ be friendly at level $j=\flev(u,x_{p-1})<p-1$.  Since there is no edge between $\{x_{p-1},x_p,x_{p+1}\}$ and $\{u^+,u^-\}$, we can apply \cref{friendly configuration} to $P_u$ and $P_x=xQ_xx_{p+1}x_px_{p-1}P_{x_{p-1}}$. Thus $j+1\in \{1,p,2p+1\}$. $j=0$ and $2p+p-2j-1=3p-1\in \{p,2p+1\}$, a contradiction.     
\end{proof}
Therefore $G$ contains $T_{n-2}^p$ as an induced subgraph. Since $G$ is deficiency-critical and $\kd(G)=\kd(T_{n-2}^p)=3n-8$, $G=T_{n-2}^p$. Hence $T_{n-2}^p$ is the only deficiency-critical graph attaining this extremum. 
\end{proof}

\onlyoneevenboneevenbone*
\begin{proof}
(1) Recall that $N\in A=\{2p\}$. Thus $N=2p$. 
\begin{claim}\label{Z_2p}
$|Z_{2p}|\le (m-2)(n-3)$.		
\end{claim}	
\begin{proof}
Suppose not, then $s=|X_{2p-1}|>m-2$ according to \cref{|X_i-1|}. Let $P_x$ and $P_{x^2}$ be friendly at level $\flev(x,x^2)$. then there exists an edge $x_jx^2_j$ for some $j>\flev(x,x^2)$, otherwise, $V(xP_x x_{\flev(x,x^2)} P_{x^2}x^2)\cup \{x^+,x^-,x^{2+},x^{2-}\}$ forms an odd bone. Let $k=\max\{j>\flev(x,x^2): x_jx^2_j\in E(G)\}$. Thus $V(xP_x x_kx^2_k P_{x^2}x^2)\cup \{x^+,x^-,x^{2+},x^{2-}\}$ forms the bone $B_{4p-2k}$. Thus $k=p>\flev(x,x^2)$. By symmetry, $\flev(x^i,x^j)<p$ for any distinct $i,j\in [s]$. 
For distinct $i,j\in [s]$, 
let $R_{ij}$ be a shortest path connecting $x^i_p$ and $x^j_p$ with all internal vertices in $L_0\cup \cdots \cup L_{p-1}$. Since  $p>\flev(x^i,x^j)$, $x^i_{\ell}x^j_{\ell+1},x^i_{\ell+1}x^j_{\ell}\not\in E(G)$ for any $\ell\ge p$. Thus $V(x^iP_{x^i}x^i_p R_{ij}x^j_p P_{x^j}x^j)\cup \{x^{i+},x^{i-},x^{j+},x^{j-}\}$ forms the bone $B_{2p+r_{ij}}$, where $r_{ij}$ is the number of internal vertices of $R_{ij}$. It forces that $r_{ij}=0$ and $x^i_px^j_p\in E(G)$. Therefore, $\{x^i_p:i\in [s]\}$ forms a clique of size $s$. For any $i\in [2,s]$, $x_{p-1}x^i_p\in E(G)$, otherwise $V(x_{p}P_xx_0)\cup \{x_0^+,x_0^-,x_{p+1},x^i_p\}$ forms the bone $B_{p+1}$, a contradiction. Therefore, $\{x_{p-1},x^i_p:i\in [s]\}$ forms a clique of size $s+1\ge m$, contradicts the assumption  that $G$ is $K_m$-free.
\end{proof}
According to \cref{Z_1}, \cref{Z_2p} and \cref{Z_inot=0}, we have
$$\kd(G)=|Z_1|+|Z_{2p}|\le (m-1)(n-3)+1.$$

(2)
Now $N\in A$ is even. Recall that $x\in X_{N-1}$ is a snail head.
\begin{claim}\label{N(x_i)}
For any selected path $P_x=xx_{N-2}x_{N-3}\dots x_1x_0$, we have $N(x_2)\subseteq\{x_3\}\cup N(x_0)$ and $N(x_{N-3})\subseteq \{x_{N-4}\}\cup N(x)$.
Additionally, For any $i\in \{1,3,4,\cdots, N-4, N-2\}$, $N(x_i)=\{x_{i-1},x_{i+1}\}$.
\end{claim}
\begin{proof}
For $i=2$ or $N-3$, by symmetry, suppose to the contrary that there exists a vertex $v\in N(x_2)$ but  $v\not \in \{x_3\}\cup N(x_0)$. Then $N(v)\cap \{x_1,x_3\}=\emptyset$ as $G$ is $K_3$-free. Thus $\{x_0,x_1,x_2,v,x_3,x_0^+,x_0^-\}$ forms the bone $B_3$, a contradiction.	
For $i\in \{1,3,4,\cdots, N-4, N-2\}$, suppose to the contrary that there exists a vertex $v\not\in \{x_{i-1},x_{i+1}\}$ adjacent to $x_i$. Thus $vx_{i-1},vx_{i+1}\notin E(G)$ as $G$ is $K_3$-free. 
If $vx_{i-2}\notin E(G)$, then $V(x_iP_xx_0)\cup \{x_0^+,x_0^-,v,x_{i+1}\}$ forms the bone $B_{i+1}$, and, otherwise, $V(x_{i-2}P_xx_0)\cup \{x_0^+,x_0^-,v,x_{i-1}\}$ forms the bone $B_{i-1}$. In a word, $i$ is odd. 
Similarly, if $vx_{i+2}\notin E(G)$, then $V(xP_xx_i)\cup \{x^+,x^-,v,x_{i-1}\}$ forms the bone $B_{N-i}$, and, otherwise, $V(xP_xx_{i+2})\cup \{x^+,x^-,v,x_{i+1}\}$ forms the bone $B_{N-i-2}$. Thus $N-i$ is even. This leads to a contradiction as $N=i+(N-i)$ is even.
\end{proof}

Recall that for any $v\in L_i$, $C(v)=N(v)\cap L_{i+1}$. 
Since $G$ is $K_3$-free, $C(v)$ is a stable set.
For even $i>0$, $|C(v)|\le 1$;  otherwise, $V(vP_vx_0)\cup C(v)\cup \{x_0^+,x_0^-\}$ forms an odd bone $B_{i+1}$, a contradiction.
For any $u\in C(x_{N-3})$, $ux_{N-1}\in E(G)$. Applying \cref{N(x_i)} to the updated path $P'_x=xux_{N-3}x_{N-4}\cdots x_0$, we have $N(u)=\{x_{N-3},x\}$. 
\begin{claim}\label{N(x_0)}
For any $u\in C(x_0)$, $C(u)\subseteq \{x_2\}$. 	
\end{claim}
\begin{proof}
If $ux_2\in E(G)$, then by applying \cref{N(x_i)} to the updated path $P'_x=xP_x x_2ux_0$, we have $N(u)=\{x_0,x_2\}$ and thus $C(u)=\{x_2\}$. The claim is true.
Now suppose to the contrary that there exists a vertex $u\in C(x_0)$ non-adjacent to $x_2$ such that $C(u)\neq \emptyset$.
Starting from $u_1:=u$, choose $u_i$ for $i\ge 2$ successively such that $u_i\in C(u_{i-1})$. The process terminates when $x_{i+1}\in C(u_i)$ or $|C(u_i)|\neq 1$. 
Let $u_p$ be the last chosen vertex. 
If $C(u_p)=\emptyset$, then $p\ge 2$ and $C(u_{p-1})=\{u_p\}$. Thus $G-\{u_p,u_{p-1}\}$ is still connected as each vertex can be connected to $x_0$ by the levelling. Since $G$ ha s a matching $M=M'\cup \{u_p,u_{p-1}\}$ where $M'$ is a maximum matching of $G-\{u_p,u_{p-1}\}$, we have $\kd(G-\{u_p,u_{p-1})\ge \kd(G)$, which contradicts to that $G$ is deficiency-critical.

Now assume $C(u_p)\not=\emptyset$. If $x_{p+1}\in C(u_p)$, then $p\neq 1$ as $ux_2\not \in E(G)$.  This forces that $x_{p+1}=x_{N-1}$ according to \cref{N(x_i)}. Thus $p=N-2$ is even. Consequently, $C(u_p)=\{x\}$ as $x\in C(u_p)$ and $|C(u_p)|\le 1$. This yields that $G-\{u_p,u_{p-1}\}$ is still connected, a contradiction.

Now we may assume that $x_{p+1}\notin C(u_p)$ and $|C(u_p)|\ge 2$. Then $p$ must be odd as $u_p\in L_p$. 
We claim that $c_p:=|C(u_p)-N(x)|\ge 2$. We may assume $C(u_p)\cap N(x)\neq \emptyset$ as $|C(u_p)|\ge 2$. Thus $p+1\in \{N-2,N-1,N\}$ as $C(u_p)\subseteq L_{p+1}$ and $N(x)\subseteq L_{N-2}\cup L_{N-1}\cup L_N$. 
By the parity of $p$ and $N$, we have $p+1=N$ or $N-2$.

If $c_p=1$, then there exist $v\in C(u_p)\cap N(x)$ and $w\in C(u_p)\setminus N(x)$ as $|C(u_p)|\ge 2$. $vw\notin E(G)$ as $C(u_p)$ is stable. Since $u_p$ exists, $xu_{p-1}\notin E(G)$. Thus $\{x,v,u_p,x^+,x^-,w,u_{p-1}\}$ forms the bone $B_3$, a contradiction. 
If $c_p=0$, take arbitrarily $v\in C(u_p)$. Thus $vx\in E(G)$ as $|C(u_p)-N(x)|=0$. If $p+1=N$, then $C(v)\subseteq L_{N+1}=\emptyset$. If $p+1=N-2$, then $C(v)=\{x\}$ as $x\in C(v)$ and $|C(v)|\le 1$. In a word, $C(v)\subseteq \{x\}$. 
Hence $G-\{u_p,v\}$ is connected since $C(v)\subseteq \{x\}$ and each vertex in $C(u_p)-\{v\}$ is adjacent to $x$. This leads to a contradiction as $G$ is deficiency-critical.

Therefore, $c_p\ge 2$ and there exists $u_p^+,u_p^-\in C(u_p)-N(x)$. 
According to the information of $N(x_i)$ by \cref{N(x_i)}, $V(u_pu_{p-1}\cdots u_2 u x_0 P_x x)\cup \{x^+,x^-,u_p^+,u_p^-\}$ forms an odd bone $B_{N+p}$, a contradiction.
\end{proof}

According to Claims~\ref{N(x_i)} and \ref{N(x_0)}, we can partition $V(G)$ into three sets: $V_1=N[x_0]-\{x_1\}$, $V_2=N[x_{N-1}]-\{x_{N-2}\}$, and $V_3=\{x_1,x_2,\cdots, x_{N-2}\}$. For $i=1,2$,  we observe that $\kd(G[V_i])\le n-3$ as $G[V_i]$ is an induced star with at most $n-2$ leaves. Additionally, $\kd(G[V_3])=0$ as $G[V_3]$ is an induced path with $N-2$ (an even number) vertices. Consequently, we have 
\[\kd(G)\le \kd(G[V_1])+\kd(G[V_2])+\kd(G[V_3])\le 2n-6.\]	
\end{proof}

\section{Remarks and Discussions}
Notably, all the extremums in Theorems \ref{main}, \ref{two odd bone} and \ref{only one even bone} (1) are congruent to $1$ modulo $(n-3)$, which naturally motivates the following problem.
\begin{prop}
Given integers $m,n>3$ and set $A\subseteq \{2,3,4,5\cdots\}$, let $\mathcal G(A;m,n)$ be the class of connected $(K_m,K_{1,n})$-free graphs in which every bone is of the form $B_i$ for some $i\in A$.
Suppose  that
\[d:=\max\{\kd(G):G\in \mathcal G(A;m,n)\}\] is finite. Is it true that  $d\equiv 1 \pmod {n-3}$? 	
\end{prop}

\noindent{\bf Acknowledgements}\\
This work was supported by the National Key Research and Development Program of China
(2023YFA1010203), the National Natural Science Foundation of China (No. 12471336), and
the Innovation Program for Quantum Science and Technology (2021ZD0302902).

\noindent{\bf Declaration of competing interest}\\
No conﬂﬂict of interest exits in the submission of this manuscript, and manuscript is approved by all authors for publication.
 I would like to declare on behalf of my co-authors that the work described was original research that has not been 
published previously, and not under consideration for publication elsewhere, in whole or in part.

\noindent{\bf Data availability}\\
No data was used for the research described in the article.


\begin{thebibliography}{99}
\bibitem{abcmmrw}
T. Abrishami, M. Briański, J. Czyżewska, R. McCarty, M. Milanič, P. Rzążewski, and B. Walczak, Excluding a clique or a biclique in graphs of bounded induced matching treewidth, arXiv:2405.04617.

\bibitem{crst}
M. Chudnovsky, N. Robertson, P. Seymour and R. Thomas, $K_4$-free graphs with no odd holes, J. Combin. Theory Ser. B. {\bf 100} (2010) 313--331.

\bibitem{cs}
M. Chudnovsky and P. Seymour, Proof of a conjecture of Plummer and Zha, J. Graph Theory. {\bf 103} (2023) 437--450.

\bibitem{css}
M. Chudnovsky, A. Scott and P. Seymour, Induced subgraphs of graphs with large chromatic number. III. Long holes, Combinatorica {\bf37} (2017) 1057-1072.

\bibitem{csss}
M. Chudnovsky, A. Scott, P. Seymour and S. Spirkl, Induced subgraphs of graphs with large chromatic number. VIII. Long odd holes. J. Combin. Theory Ser. B {\bf 140} (2020) 84-97.


\bibitem{dkkmmsw}
C. Dallard, M. Krnc, O. Kwon, M. Milani\v c, A. Munaro, K. \v Storgel and S. Wiederrecht, Treewidth versus clique number. IV. Tree-independence number of graphs excluding an induced star, arXiv:2402.11222.

\bibitem{dm}
X. Du and R. McCarty, a survey of degree-boundedness, arXiv:2403.05737v2.

\bibitem{dms}
C. Dallard, M. Milani\v c{} and K. \v Storgel, Treewidth versus clique number. I. Graph classes with a forbidden structure, SIAM J. Discrete Math. {\bf 35} (2021), no.~4, 2618--2646.

\bibitem{fklops}
S. Fujita, K. Kawarabayashi, C. L. Lucchesi, K. Ota, M. Plummer and A. Saito. A pair of forbidden subgraphs and perfect matchings. J. Combin. Theory Ser. B {\bf 96}, (2006) 315--324.

\bibitem{g}
A. Gy\'arf\'as, Problems from the world surrounding perfect graphs, in: Proceedings of the Inter- national Conference on Combinatorial Analysis and its Applications, Pokrzywna, 1985, Zastos. Mat.
{\bf 19} (1987) 413--441.

\bibitem{gz}
A. Gy\'arf\'as and M. Zaker, On $(\delta,\chi)$-bounded families of graphs, Electron. J. Combin. {\bf 18} (2011), no.~1, Paper 108, 8 pp. 

\bibitem{jpr}
M. J\"unger, W. R. Pulleyblank and G. Reinelt. On partitioning the edges of graphs into connected subgraphs, J. Graph Theory {\bf 9} (1985) 539-549.

\bibitem{l}
M. Las Vergnas. A note on matchings in graphs, Colloque sur la Th\'eorie des Graphes (Paris 1974), Cahiers Centre \'Etudes Rech. Op\'er. {\bf 17} (1975) 257-260.

\bibitem {s}
D. P. Sumner. $1$-factors and antifactor sets, J. London Math. Soc. {\bf 13} (1976) 351-359.

\bibitem{sr}
I. Schiermeyer and B. Randerath, Polynomial $\chi$-binding functions and forbidden induced subgraphs: a survey, Graphs Combin. {\bf 35} (2019), no.~1, 1--31.

\bibitem{ss}
A. Scott and  P. Seymour, Induced subgraphs of graphs with large chromatic number. I. Odd holes, J. Combin. Theory Ser. B {\bf 121} (2016) 68-84.

\bibitem{ss20}
A. Scott and P. Seymour, A survey of $\chi$-boundedness, J. Graph Theory {\bf 95} (2020), no.~3, 473--504.

\bibitem{wxx}
D. Wu, B.~G. Xu and Y. Xu, The chromatic number of heptagraphs, J. Graph Theory {\bf 106} (2024), no.~3, 711--736.
\end{thebibliography}
\end{document}